\documentclass[letterpaper,11pt]{amsart}

\usepackage{latexsym}

\RequirePackage[OT1]{fontenc}
\RequirePackage{amsmath}
\RequirePackage[colorlinks,citecolor=blue,urlcolor=blue]{hyperref}
\usepackage{exscale}
\usepackage{amssymb}
\usepackage{amscd}
\usepackage{bbm}
\usepackage{color}
\usepackage{graphicx}
\usepackage{mathrsfs}
\usepackage{comment}
\usepackage{multicol}

\usepackage{enumitem} 

\DeclareGraphicsExtensions{.eps, .ps}

\newcommand{\Leb}{{\rm Leb}}

\renewcommand{\to}{\rightarrow}
\newcommand{\upto}{\uparrow}

\renewcommand{\H}{\mathcal{H}}
\newcommand{\N}{\mathbb{N}}

\newcommand{\bw}{\wedge}

\def\build#1_#2^#3{\mathrel{\mathop{\kern 0pt#1}\limits_{#2}^{#3}}}

\newcommand{\bN}{\mathbb{N}}

\newcommand{\bZ}{\mathbb{Z}}

\newcommand{\cG}{\mathcal{G}}

\newcommand{\cI}{\mathcal{I}}

\newcommand{\cL}{\mathcal{L}}

\newcommand{\cT}{\mathcal{T}}

\newcommand{\whH}{\widehat{\mathcal{H}}}

\newcommand{\tT}{\mathtt{T}}
\newcommand{\tts}{\mathtt{s}}

\def\cf{\mathbf{1}}
\def\Be{\mathbf{e}}

\def\Pr{\mathbb{P}}
\def\EV{\mathbb{E}}

\newcommand{\T}{\mathcal{T}}
\def\scH{\mathscr{H}}

\def\tail{\textnormal{tail}}
\def\cl{\textnormal{cl}}

\def\EIG{\mathtt{EIG}}

\newcommand{\Restrict}[2]{#1\big|_{#2}}
\newcommand{\restrict}[2]{#1|_{#2}}
\newcommand{\fringe}[2]{F_{#2}(#1)}

\theoremstyle{plain}

\newtheorem{theorem}{Theorem}

\newtheorem{proposition}{Proposition}
\newtheorem{corollary}{Corollary}
\newtheorem{lemma}{Lemma}
\newtheorem{conjecture}{Conjecture}

\theoremstyle{definition}
\newtheorem{definition}{Definition}
\newtheorem{example}{Example}

\theoremstyle{remark}

\title[Exchangeable hierarchies]{A representation of exchangeable hierarchies by sampling from random real trees}

\author[N.~Forman]{Noah Forman$^1$}
\address{$^1$ Department of Mathematics\\ University of Washington\\ Seattle, WA 98195}
\author[C.~Haulk]{Chris Haulk$^2$}
\address{$^2$ Google}
\author[J.~Pitman]{Jim Pitman$^3$}
\address{$^3$ Statistics Department\\ Evans Hall\\ University of California, Berkeley\\ Berkeley, CA 94720}

\date{\today}

\thanks{Research supported in part by NSF grants DMS-0806118 and DMS-1444084 and EPSRC grant EP/K029797/1}

\subjclass[2010]{Primary 60G09; Secondary 60C05, 62B05}
\keywords{exchangeable, hierarchy, total partition, random composition, random partition, continuum random tree, weighted real tree}

\begin{document}
\begin{abstract}
A {\em hierarchy } on a set $S$, also called a {\em total partition of $S$}, is a collection $\H$ of subsets of $S$ such that $S \in \H$, each singleton subset of $S$ belongs to $\H$, and if $A, B \in \H$ then $A \cap B$ equals either $A$ or $B$ or $\varnothing$.  Every exchangeable random hierarchy of positive integers has the same distribution as a random hierarchy $\H$ associated as follows with a random real tree $\T$ equipped with root element $0$ and a random probability distribution $p$ on the 
Borel subsets of $\T$: given $(\T,p)$, let $t_1,t_2, \ldots$ be independent and identically distributed according to $p$, and let $\H$ comprise all singleton subsets of $\N$, and every subset of the form $\{j\colon t_j \in F(x)\}$ as $x$ ranges over $\T$, where $F(x)$ is the fringe subtree of $\T$ rooted at $x$. There is also the alternative characterization: every exchangeable random hierarchy of positive integers has the same distribution as a random hierarchy $\H$ derived as follows from a random hierarchy $\mathscr{H}$ on $[0,1]$ and a family $(U_j)$ of i.i.d.\ Uniform[0,1] random variables independent of $\mathscr{H}$: let $\H$ comprise all sets of the form $\{j\colon U_j \in B\}$ as $B$ ranges over the members of $\mathscr{H}$. 
\end{abstract}

\ \vspace{-38pt}

\maketitle

\section{Introduction}\label{sect-intro}

\begin{definition}\label{defn-hierarchy}
 A {\em hierarchy} on a finite set $S$ is a set $\H$ of subsets of $S$ for which 
 \begin{enumerate}[label=(\alph*), ref=(\alph*)]
  \item if $A, B \in \H$ then $A \cap B$ equals either $A$ or $B$ or $\varnothing$, and
  \item $S \in \H$, $\{s\} \in \H$ for all $s \in S$, and $\varnothing \in \H$.  
 \end{enumerate}
\end{definition}

Hierarchies are known by several other names, including {\em total partitions} and \emph{laminar families}.  For brevity we use the term \emph{hierarchy} throughout the paper. Less formally, a hierarchy describes a scheme for recursively partitioning a set $S$ into
finer and finer subsets, down to singletons. 
 Alternatively, a hierarchy describes a process of coalescence, wherein the singleton subsets of $S$ recursively coagulate to reconstitute the set $S$.  We emphasize that \emph{time} plays no role in our definition of a hierarchy: a hierarchy encodes the \emph{contents} of the blocks of some process of fragmentation (or coagulation), but includes no information about the order in which these blocks appear.

Hierarchies on $[n]$ correspond to certain trees.  If $\mathtt{T}$ is a tree 
\begin{itemize}
 \item{with $n$ leaves, labeled by distinct elements of $[n]$;}
 \item{having a distinguished vertex called the root, which is not a leaf;}
 \item{with no internal vertices of degree two, except possibly the root;}
 \item{and no edge lengths or planar embedding}
\end{itemize}
then the map
\[
 \mathtt{T} \mapsto \big\{\{j \in [n]\colon v \mbox{ on path from leaf }j \mbox{ to root}\}\colon v \in V(\mathtt{T})\big\}\cup \{\varnothing\}
\]
sends $\mathtt{T}$ to a hierarchy on $[n]$.  Here, $V(\mathtt{T})$ denotes the set of vertices  of $\mathtt{T}$ (including root and leaves).  This map is a bijection, and we say that $\mathtt{T}$ is the \emph{graph} of the corresponding hierarchy.

For a hierarchy $\H$ on a set $S$ and $S_0 \subseteq S$, the \emph{restriction of $\H$ to $S_0$} is
\begin{equation}\label{restriction}
 \Restrict{\H}{ S_0} := \{H \cap S_0\colon H \in \H\},
\end{equation}
which is a hierarchy on $S_0$. We equip the set of hierarchies on a finite set $S$ with the discrete $\sigma$-algebra. 

\begin{figure}[htbp] 
   \centering
   \includegraphics[width=4in]{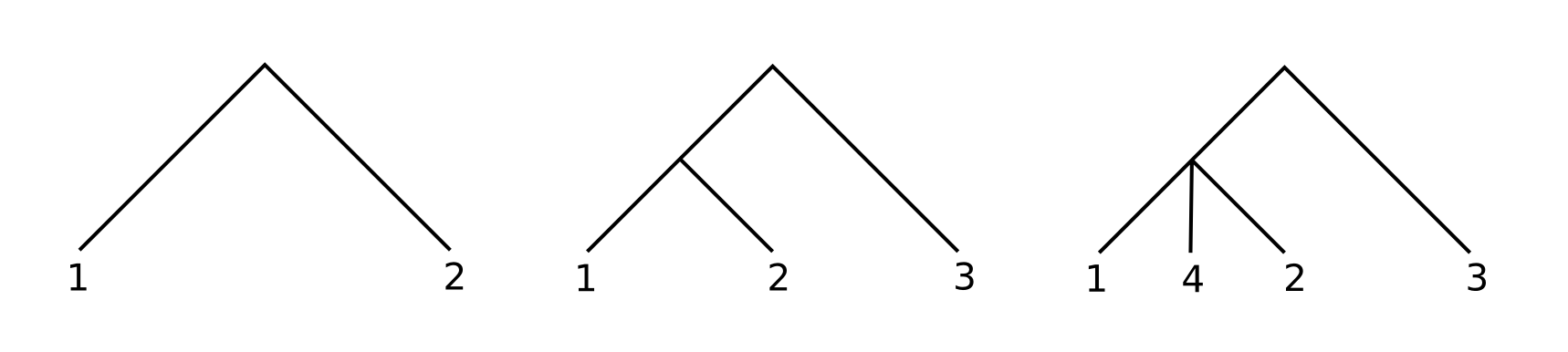} 
   \caption{Right: the graph of $\H = \{\{1,2,4\}\} \cup\Xi([4])$, where $\Xi$ is the trivial hierarchy of \eqref{defn-xi}. Center: the graph of $\Restrict{\H}{[3]}$. Left: the graph of $\Restrict{\H}{[2]}$.
   \label{fig-conistency}}
\end{figure}


\begin{definition}\label{def:N_hier}
 A \emph{hierarchy on $\bN$} is a sequence $(\H_n,\,n\ge1)$, with each $\H_n$ a hierarchy on $[n]$, with the consistency condition that $\H_n = \Restrict{\H_{n+1}}{[n]}$ for each $n\ge 1$. We take the $\sigma$-algebra on the space of such hierarchies to be that generated by the projection maps $(\H_n,\,n\ge1) \mapsto \H_k$ for $k\in\bN$.
\end{definition}




Permutations act on hierarchies by relabeling the contents of constituent sets: if $\H$ is a hierarchy on $[n]$ and $\sigma$ a permutation of $[n]$, then 
\[
 \sigma( \H ):= \big\{ \{\sigma(j)\colon j \in H\}\colon H \in \H\big\}.
\]
An {\em exchangeable hierarchy on $\mathbb{N}$} is a random hierarchy $(\H_n,\, n \geq 1)$ on $\mathbb{N}$ for which, for every $n\in\N$ and every permutation $\sigma$ of $[n]$, there is the distributional equality
\begin{equation}
 \sigma(\H_n) \stackrel{d}{=} \H_n.
\end{equation}

Recall that a sequence of random variables $(X_n,\,n\ge1)$ is exchangeable if $(X_n,\,n \geq 1) \stackrel{d}{=} (X_{\sigma(n)},\,n\ge1)$ for every finite permutation $\sigma$ of $\N$. Let $\tail(X_n)$ denote the associated tail $\sigma$-algebra.
\begin{theorem}[de Finetti's Theorem \cite{MR2161313}] \label{theorem:de-finetti}
 If $(X_n,\, n\ge1)$ is an exchangeable sequence of random elements of a Borel measurable space $(S,\mathcal{S})$, then there is a random probability measure $\nu$ on $(S,\mathcal{S})$, measurable with respect to $\tail(X_n)$, 
 for which 
 $\nu^{\bN}$ is a regular conditional distribution (r.c.d.) for $(X_n,\, n\ge1)$ on $\tail(X_n)$. I.e.\ conditionally given $\tail(X_n)$, the sequence $(X_n,\, n\ge1)$ is i.i.d.\ with law $\nu$.
\end{theorem}

This random measure $\nu$ is called the \emph{driving measure} of the sequence $(X_n)$. For a hierarchies version of this theorem we require corresponding notions of independence and the tail $\sigma$-algebra.

\begin{definition}\label{def:EIG}
 A hierarchy $\H$ on a set $S$ is \emph{independently generated} if for each vector $(A_1,\ldots ,A_k)$ of disjoint, finite subsets of $S$, the restrictions $\big(\restrict{\H}{A_1},...,\restrict{\H}{A_k}\big)$ of $\H$ to these subsets are mutually independent. We write \emph{e.i.g.}\ to abbreviate ``exchangeable and independently generated.'' Denote by $\EIG$ the set of distributions of e.i.g.\ hierarchies on $\N$.
\end{definition}

\begin{definition}\label{def:tail}
 The \emph{tail $\sigma$-algebra} of a random hierarchy $(\H_n,\,n\ge1)$ on $\N$ is
 \[
  \tail(\H_n) := \bigcap_{n\in\N} \sigma\big(\restrict{\H}{\{n\}},\restrict{\H}{\{n,n+1\}}, \restrict{\H}{\{n, n+1, n+2\}},\ldots\big).
 \]
\end{definition}

We offer the following hierarchies analogue to de Finetti's Theorem.

\begin{theorem}\label{thm:EIG}
 \begin{enumerate}[label=(\roman*), ref=(\roman*)]
  \item $\EIG$ is the set of extreme points in the convex set of probability distributions of exchangeable hierarchies on $\N$. \label{item:EIG:convex}
  \item  If $(\H_n,\,n\ge1)$ is an exchangeable hierarchy on $\N$ then, given $\tail(\H_n)$, $\H$ is conditionally e.i.g.. That is, there is a random law $\cL\in\EIG$ that is a r.c.d.\ for $(\H_n)$ given $\tail(\H_n)$. \label{item:EIG:rcd}
  \item If $\mathcal{M}$ is an $\EIG$-valued r.c.d.\ for $(\H_n,\,n\ge1)$ given some $\sigma$-algebra, then $\mathcal{M} = \cL$ almost surely. \label{item:EIG:unique}
 \end{enumerate}
\end{theorem}


This theorem is proved in Section \ref{sec:EIG_proof}.


\subsection{Kingman-type descriptions}

In preparation for more concrete descriptions of exchangeable hierarchies, we state a version of Kingman's representation theorem for exchangeable partitions.  
A random partition $\Pi$ of $\mathbb{N}$ is said to be {\em exchangeable} if the random array 
\[
 \mathbf{p}(i,j) := \cf\{i\text{ and }j\text{ are in same block of }\Pi\},\qquad i,j\in\N
\]
is exchangeable, meaning that for every $n \geq 1$ and permutation $\sigma$ of $[n]$,
\begin{equation}\label{exch-partition-eqn}
 \big( \mathbf{p}(\sigma(i), \sigma(j)),\ i,j \in [n]\big) \stackrel{d}{=} \big(\mathbf{p}(i, j),\ i,j \in [n]\big) .
\end{equation}
We can define a tail $\sigma$-algebra associated with $\Pi$ in a manner analogous to Definition \ref{def:tail}.

\begin{theorem}[Kingman's representation theorem \cite{MR509954}]\label{thm:partn_Kingman}
 Suppose that $\Pi$ is an exchangeable partition of $\bN$ and that the probability space supports a sequence $(s_i,\,i\ge 1)$ of i.i.d.\ Uniform$[0,1]$ random variables independent of $\Pi$.  Then there is a $\tail(\Pi)$-measurable random open subset $\mathscr{U}$ of $[0,1]$ such that for the partition $\Pi'$ defined by
 \[
  \{i\text{ and }j\text{ in same block of }\Pi'\} =  \{s_i\text{ and }s_j\text{ are in same interval in }\mathscr{U}\},
 \]  
 we get equality of joint distributions $(\Pi, \mathscr{U}) \stackrel{d}{=} (\Pi', \mathscr{U})$.
\end{theorem}

In this setting, we view the set of maximal open subintervals of $\mathscr{U}$ as forming an \emph{interval partition} of $[0,1]$, though this partition may be incomplete in the sense that the complement $[0,1]\setminus \mathscr{U}$ may have positive Lebesgue measure. See \cite[Section 17]{MR883646} or \cite[Chapter 4]{MR2245368}.


\begin{definition}\label{def:derived_by_sampling}
 \begin{enumerate}[label=(\roman*),ref=(\roman*)]
  \item If $S$ is an infinite set and $\scH$ satisfies properties (a) and (b) of Definition \ref{defn-hierarchy}, then we call $\scH$ an \emph{$\infty$-hierarchy on $S$}. We will not refer to $\infty$-hierarchies on $\bN$, in which setting we prefer Definition \ref{def:N_hier}.
  \item Suppose $\scH$ is an $\infty$-hierarchy on a set $S$, with $p$ a probability distribution on $S$, and $(s_i,\,i\ge1)$ i.i.d.\ with law $p$. Setting
 \begin{equation}
  \H_n := \big\{\{ i \in [n]\colon s_i\in A\} \colon A\in\scH\big\} \cup \big\{\{i\}\colon i\ge 1\big\} \qquad \text{for } n\ge1,
 \end{equation}
 we say $(\H_n,\,n\ge1)$ is the \emph{hierarchy derived from $\scH$ by samples $(s_i)$}, or is \emph{derived by sampling from $(\scH,p)$}. Let $\Theta(\scH,p)$ denote the law of $(\H_n)$. When $\scH$ is a deterministic hierarchy, $\Theta(\scH,p)$ is an e.i.g.\ law.\label{item:hier_sampling}
  \item For real numbers $a<b$, an \emph{interval hierarchy on $[a,b)$} is an $\infty$-hierarchy on $[a,b)$ in which each non-singleton block is an interval $[c,d)$ as well. We take the $\sigma$-algebra on the space of such hierarchies to be the least $\sigma$-algebra under which the restrictions $\scH\mapsto \Restrict{\scH}{A}$ to finite sets $A\subset [a,b)$ are measurable.\label{item:int_hier}
 \end{enumerate}
\end{definition}

\begin{figure}[thbp] 
   \centering
   \includegraphics[width=4in]{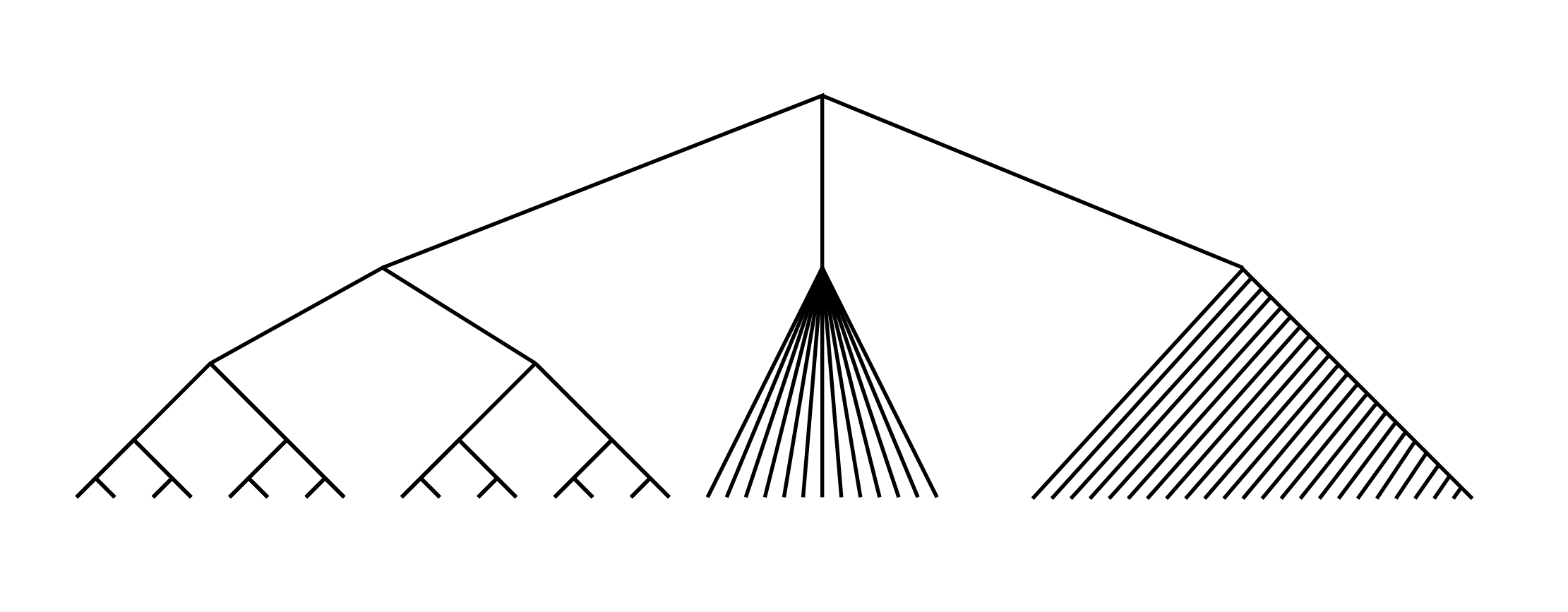} 
   \caption{Approximate graph of the hierarchy in Example \ref{eg:naive}, with leaf labels omitted.}
   \label{fig:naive}
\end{figure}

\begin{example}\label{eg:naive}
 Let $\mathscr{H}$ denote the following interval hierarchy on $[0,3)$:
 \begin{equation*}
 \begin{split}
  \mathscr{H} &:= \bigcup\nolimits_{n \geq 1} \left\{\left[\frac{j}{2^n}, \frac{j+1}{2^n}\right)\colon 0 \leq j \leq 2^n-1\right\} \cup \{[x,3)\colon 2 < x < 3\}\\
 		&\qquad \cup \{\{x\}\colon x\in [0,3)\} \cup  \{[0,1), [1,2), [2,3), [0,3),\varnothing\}.
 \end{split}
 \end{equation*}
 Let $(\H_n,\,n\ge 1)$ be derived by sampling from $(\scH,\text{Uniform}[0,3))$. This is illustrated in Figure \ref{fig:naive}.
\end{example}


To help build intuition around hierarchies, we mention the following.

\noindent\textbf{Na\"{\i}ve conjecture.} 
 The three phenomena exhibited in Figure \ref{fig:naive} -- iterative splitting into non-singleton blocks, broom-like explosion into singletons, and comb-like erosion of singletons –- are the basic building blocks out of which every exchangeable, independently generated hierarchy is made.
 
A similar characterization arises in the theory of fragmentation processes \cite{MR1899456,MR3141802}. We will not prove this conjecture here, but we propose a formalization of it at the end of this article.

Here is our first Kingman-type result for hierarchies. Let $\Leb$ denote Lebesgue measure.

\begin{theorem}\label{thm:interval_Kingman}
 If $(\H_n,\,n\ge1)$ is an exchangeable random hierarchy on $\bN$ then there is an $(\H_n)$-measurable random interval hierarchy $\mathscr{H}$ on $[0,1)$ such that $\Theta(\scH,\Leb)$ is a r.c.d.\ for $(\H_n)$ on $\tail(\H_n)$. Moreover,
 $$\EIG = \{\Theta(\scH,\Leb)\colon \scH\text{ an interval hierarchy on }[0,1)\}.$$
\end{theorem}

We offer a second Kingman-type description in terms of sampling from real trees. 
%
For now, we state only the definitions needed for our main result. Examples and a general construction may be found in Section \ref{sec:line_break}.

\begin{definition} \label{def:real_tree}
 A \emph{segment} of a metric space $X$ is the image of an isometry on a real interval, \mbox{$\alpha\colon [a, b] \mapsto X$}.  The \emph{endpoints} of the the segment are $\alpha(a)$ and $\alpha(b)$.  A \emph{real tree} is a metric space $(\T,d)$ for which 
 \begin{enumerate}[label=(\alph*), ref=(\alph*)]
  \item for every pair $x,y\in \T$ there is a unique segment with endpoints $x$ and $y$, denoted $[[x,y]]_{\T}$, or simply $[[x,y]]$ where there is no ambiguity (note that $[[x,x]]_{\T} = \{x\}$);
  \item if two segments of $\T$ intersect in a single point, and this point is an endpoint of both, then the union of these two segments is again a segment;
  \item if a segment contains distinct points $u$, $v$ then it contains $[[u,v]]_{\T}$; and
  \item if the intersection of two segments contains at least two distinct points, then this intersection is a segment.
 \end{enumerate}
 
 A \emph{rooted, weighted real tree} is a quadruple $(\T,d,r,p)$, where $(\T,d)$ is a real tree, $r\in \T$ is a distinguished vertex called the \emph{root}, and $p$ is a probability distribution on $\T$ with respect to the Borel $\sigma$-algebra generated by $d$.
 
 A point $x\in \T$ is a \emph{leaf} if it is an endpoint of every segment to which it belongs. It is a \emph{branch point} if there exists three non-trivial segments with endpoint $x$ whose pairwise intersections all equal $\{x\}$. The complement of the set of leaves is the \emph{skeleton} of the tree.
\end{definition}

Random real trees were first studied by Aldous \cite{MR1085326,MR1207226}. For more on this topic see the course notes \cite{MR2351587}, in which properties (c) and (d) are derived as consequences of (a) and (b). Also see \cite{MR2203728}. 
For $(\T,d,r)$ a rooted real tree and $x\in\T$, the \emph{fringe subtree} of $\T$ rooted at $x$ is
\begin{equation}
 \fringe{x}{\T}:= \{ y \in \T\colon x \in [[r,y]]\}.\label{fringe-subtree-defn}
\end{equation}
For the purpose of the following, we denote
\begin{equation}
 \scH(\T,d,r) := \{\fringe{x}{\T}\colon x\in\T\} \cup \{\{x\}\colon x\in\T\} \cup \{\varnothing\}.
\end{equation} 
This is easily shown to be an $\infty$-hierarchy on $\T$. 
We adopt the abbreviation $\Theta(\T,d,r,p)$ to denote $\Theta(\scH(\T,d,r),p)$. 
%

\begin{theorem}\label{thm:tree_constr}
 Let $(\H_n,\,n\ge1)$ be an exchangeable random hierarchy on $\N$. Then there exists an $(\H_n)$-measurable random rooted, weighted real tree $(\T,d,r,p)$ such that $\Theta(\T,d,r,p)$ is a r.c.d.\ for $(\H_n)$ given $\tail(\H_n)$.
\end{theorem}

Theorems \ref{thm:interval_Kingman} and \ref{thm:tree_constr} are not as complete as Kingman's description of partitions in the following sense. Given two interval partitions $\mathscr{U}$ and $\mathscr{U}'$, it is easy to tell whether both correspond to the same law on exchangeable partitions of $\N$: looked at the ranked block sizes of the partitions. But it is not easy to do the same for pairs of interval hierachies or rooted, weighted real trees. 
In this vein, forthcoming work will resolve the following conjecture.

\begin{conjecture}\label{conj:IP_isom}
 There exists a class of rooted, weighted real trees such that each exchangeable hierarchy $(\H_n,\, n\ge1)$ is represented, in the sense of Theorem \ref{thm:tree_constr}, by a random member of this class that is a.s.\ unique up to isomorphism (in a suitable sense). Moreover, the isomorphism class of this random tree is $\tail(\H_n)$-measurable.
 %
\end{conjecture}

In Section \ref{sec:EIG_proof} we prove Theorem \ref{thm:EIG}.  In Section \ref{sect-preliminaries} we discuss real trees and state definitions, results from the literature, and elementary propositions needed for Theorem \ref{thm:tree_constr}, which we prove by construction in Section \ref{sec:construction}. We use this to prove Theorem \ref{thm:interval_Kingman} in Section \ref{sec:interval_results}. Section \ref{sect-complements} offers complementary discussion and miscellaneous results.

\subsection{Related work}\label{sect-bkg}

Random hierarchies of both finite and infinite sets arise naturally in a number of applications, including stochastic models for phylogenetic trees 
\cite{MR2367031, MR1704343, MR1838600, MR1681126, MR671034, MR675968, MR1880231, MR1742154}, processes of fragmentation and coalescence 
\cite{MR1434128, MR2446321,MR1771665,MR1825153,  MR2086161,MR2253162,MR2477381, MR1990057,MR2145730,MR1797304,MR1954248,MR1765005,MR2484170,MR2164028,MR2440924,MR3141802}, and statistics and machine learning
\cite{MR2606082, YeeWhyeTeh, ghahramani2010tree, Heller-Ghahramani, Meeds}.  In these applications, the object of common interest is a rooted tree which describes evolutionary relationships (in the case of phylogenetic trees) or the manner in which an object fragments into smaller pieces (as in models of fragmentation) or some notion of class membership (as in hierarchical clustering).  Such trees may have edges equipped with lengths that measure the time between speciations or fragmentation events, or some measure of dissimilarity or distance between classes. The hierarchies that we consider correspond with trees of this type \emph{without} edge lengths.

%
%

As indicated in \cite{MR1943145}, an exchangeable hierarchy $(\H_n,\,n\ge1)$ is generated by each of Bertoin's homogeneous fragmentation processes. Moreover, associated with each such fragmentation process there is a one-parameter family of self-similar fragmentations, each obtained from the homogeneous fragmentation by a suitable family of random time changes, and each generating the same random hierarchy $(\H_n,\,n\ge1)$.  An attractive feature of the self-similar fragmentations of index $\alpha <0$ is that each sample path of such a fragmentation is associated with a compact real tree \cite{MR2041829}.  The sample paths of Kingman's coalescent \cite{MR671034} can likewise be naturally identified with a compact real tree \cite{MR1765005}.  There has been considerable interest in describing real tree limits of discrete trees with edge-lengths \cite{MR1085326, MR1166406, MR1207226, MR3112436, MR1797304, MR1964956, MR3601650, MR2440924, MR2546748, MR2561439}, and Theorem \ref{thm:tree_constr} of this paper is in a similar vein. 

This work forms part of a growing list of characterizations of infinite exchangeable combinatorial objects by de Finetti-type theorems.  For example, Kingman characterized exchangeable partitions of $\mathbb{N}$ \cite{MR509954}, Donelly and Joyce \cite{MR1457625} and Gnedin \cite{MR1104078} characterized composition structures, Janson characterized exchangeable posets \cite{Janson}, and Hirth characterized exchangeable ordered trees \cite{MR1982030}, 
about which we say a few words below. 
Many related de Finetti-type theorems are known \cite{MR637937, MR556418,MR577313, MR2463439,MR642729, MR2529787, MR1363235, MR2327835, MR1031426,MR1175277}, and there are excellent treatments in \cite{MR883646, MR2161313} of related material.  Such de Finetti-type results are often proved via reverse martingale convergence arguments, similar in spirit to the modern approach to de Finetti's Theorem in \cite[Chapter 4]{durrett}.  Alternate approaches use harmonic analysis \cite{MR1847948, MR1785522, MR1902597, MR2463440}, isometries of $L^2$ \cite{MR2426176}, or Choquet theory \cite{MR0076206}. The results of this paper are proved using a different approach, the key idea of which is to encode an exchangeable hierarchy using a binary array, show that this array inherits exchangeability from the hierarchy, and apply well-known characterization theorems for arrays.  A similar approach was first used by Aldous, who simplified Kingman's proof characterizing of exchangeable partitions of $\mathbb{N}$ by encoding such partitions as exchangeable sequences of real random variables \cite{MR883646}. 

In \cite[Theorem 3]{MR1207226} it is shown that if $(\mathcal{R}(k),\, k \geq 1)$ is a consistent family of exchangeable trees with edge lengths that is \emph{leaf-tight} then $(\mathcal{R}(k),\, k \geq 1)$ is derived as if by sampling from a random real tree.  (Aldous also assumes that his trees are binary, but this assumption is not essential to his proof.)  Since a hierarchy on $\mathbb{N}$ corresponds to a sequence of consistent trees \emph{without} edge lengths, the main results of this paper can be seen as a variation on this result of Aldous, showing that leaf-tightness (and indeed any pre-defined notion of distance) is not needed to obtain a de Finetti-type theorem for trees with exchangeable leaves.

In \cite{MR1990057}, it is shown that every \emph{exchangeable $\mathcal{P}$-coalescent process} corresponds to a unique \emph{flow of bridges}.  An exchangeable $\mathcal{P}-$coalescent process is a Markov process $(\Pi_t,\, t \geq 0)$ whose state space $\mathcal{P}$ is the set of partitions of $\mathbb{N}$, for which $\Pi_t$ is an exchangeable partition of $\mathbb{N}$ for every $t \geq 0$ whose increments are independent and stationary, if the notion of ``increments'' of a $\mathcal{P}$-valued function is properly understood.  This provides a de Finetti-type characterization of exchangeable coalescents.  One may ``forget'' time by setting $\H:=\{B\subset \mathbb{N}\colon B \in \Pi_t \text{ for some } t > 0\} \cup \{\mathbb{N} \}$ and thereby obtain an exchangeable hierarchy $\H$ on $\mathbb{N}$ (the notation $B \in \Pi_t$ means that $B$ is a block in the partition $\Pi_t$).  The results of Bertoin and Le Gall in \cite{MR1990057} therefore provide a de Finetti-type characterization of hierarchies that arise in this manner from exchangeable coalescents.  Due to the stationary, independent increments property, this class of hierarchies is far from including every exchangeable hierarchy, so the present work may be seen as extending the results of Bertoin and Le Gall.    

Haas and Miermont  \cite{MR2041829} provide a de Finetti-type representation of self-similar fragmentations of index $\alpha < 0$ that have no erosion or sudden loss of mass in terms of continuum trees $(\T, p)$ as follows: every such fragmentation $(F(t),\, t \geq 0)$ is derived as if from a continuum tree $(T, p)$ by setting $F(t)$ equal to the decreasing sequence of masses of connected components of $\{v \in \T\colon \textrm{ht}(v) >t\}$ where ht($v$) denotes the distance from $v$ to the root of $\T$.  This is proved by introducing a family $(R(k),\, k \geq 1)$ of trees derived from an associated fragmentation $(\Pi_t)$ whose sequence of ranked limit frequencies equals $(F(t))$.  Distances in these trees $R(k)$ are related to times between dislocations in $F(t)$, and by using  self-similarity the \emph{leaf-tight} criterion of \cite{MR1207226} is checked.  The existence of the representing tree $(\T,p)$ is then a consequence of the aforementioned theorem of Aldous.  This provides a de Finetti-type theorem for self-similar fragmentations.  

Austin and Panchenko \cite{MR3230009} study processes indexed by full infinitary trees (or hierarchies) of finite depth, distributionally invariant under rooted tree automorphisms, similar to \cite[\S 13]{MR883646}, and prove a de Finetti-type result in this setting.

In \cite{MR1982030}, Hirth considers \emph{exchangeable ordered trees}, which in our terms are exchangeable hierarchies $\H$ on $\mathbb{N}$ in which, for every element $B \in \H$ besides $\mathbb{N}$, there is an associated nonnegative integer-valued ``birth time'' $N_B$ and ``death time'' $M_B$. There is also a partial order on such blocks $B$ that is unimportant for our purposes.  At the instant of its death, $B$ gives birth to subsets whose union is $B$.  Hirth provides a de Finetti-type characterization of exchangeable ordered trees using harmonic analysis techniques.  Our hierarchies are more general than Hirth's trees, since there is no ``discrete time'' associated to the elements of a hierarchy.  Our results may therefore be seen as an extension of Hirth's result using probabilistic techniques instead of harmonic analysis.  

Gufler \cite{GuflerCoalescent} describes exchangeable random semi-ultrametrics on $\N$. Ultrametric spaces can be embedded into the leaf sets of real trees. En route to describing growing genealogies as tree-valued processes, Gufler shows that exchangeable semi-ultrametrics are distributed as if by sampling from a random weighted real tree, as in Theorem \ref{thm:tree_constr} above. 
However, again, distances are inherent in the problem addressed in \cite[Theorem 1.2]{GuflerCoalescent}, whereas in our setting a metric must be introduced artificially.

Evans, Gr\"ubel, and Wakolbinger \cite[Section 5]{MR3601650} study limits of a stochastic binary tree growth procedure, using a similar notion of most recent common ancestors to that applied Section \ref{sec:MRCA} below. Moreover, \cite[Sections 6-7]{MR3601650} may be seen as an alternative to our construction in Section \ref{sec:construction}.

\section{Proof of Theorem \ref{thm:EIG}}\label{sec:EIG_proof}

%
 (i) We will prove that all e.i.g.\ laws are extreme in the set of laws of exchangeable hierarchies. The inverse, that laws that are not independently generated are not extreme, follows as an easy consequence of assertion (ii).
 
 Let $\mu$ be the distribution of an e.i.g.\ hierarchy $(\H_n,\,n\ge1)$ and consider distributions $\lambda_1$ and $\lambda_2$ of exchangeable hierarchies such that $\mu = \frac12(\lambda_1 + \lambda_2)$.  Fix a non-random hierarchy $h$ on $A = [n]$, let $B = \{n+1,\ldots, 2n\}$ and let $h'$ be the image of $h$ under $x \mapsto x+n$. Consider events $E := \{\H_n = h\}$ and $E' := \{\restrict{\H_{2n}}{B} = h'\}$.  Note that, since $(\H_n)$ is e.i.g.,
 \begin{equation*}
 \begin{split}
  \mu(E\cap E') &= \mu(E)\mu(E') = \mu(E)^2 = \frac14 (\lambda_1(E) + \lambda_2(E))^2,\\
  \text{while} \quad 
  \mu(E\cap E') &= \frac12 (\lambda_1(E\cap E') + \lambda_2(E\cap E'))\quad \text{as well}.
 \end{split}
 \end{equation*}
 It is a familiar fact that for an exchangeable Bernoulli sequence $(X_i,\,i\ge1)$ we have $\Pr(X_1 = 1, X_2 = 1) \ge \Pr(X_1 = 1)^2$, as follows from de Finetti's theorem.  Using appropriate indicator functions, this implies 
 $$\lambda_i(E\cap E') \ge \lambda_i(E)^2\qquad \text{for } i=1,2.$$
 It follows that 
 \begin{equation*}
 \begin{split}
  0 &=		\frac14(\lambda_1(E) + \lambda_2(E))^2  -  \frac12(\lambda_1(E\cap E') + \lambda_2(E\cap E'))\\
  	&\leq	\frac14(\lambda_1(E)\! +\! \lambda_2(E))^2  -  \frac12(\lambda_1(E)^2\! +\! \lambda_2(E)^2)
  	=		-\frac14 (\lambda_1(E)-\lambda_2(E))^2 \leq 0.
 \end{split}
 \end{equation*}
 Thus, $\lambda_1(E) = \lambda_2(E)$. Since this holds for every $h$, $\lambda_1 = \lambda_2 = \mu$, so $\mu$ is extreme.

(ii) Let $(\H_n,\,n\ge1)$ be an exchangeable hierarchy. 
 Let $\cG^j_N$ denote the restriction of $(\H_n,\,n\ge1)$ to $N(j-1)+[N]$, with indices shifted so that it becomes a hierarchy on $[N]$:
 $$\cG^j_N := \big\{ \{i\in [N] \colon i+N(j-1) \in A\} \colon A\in \H_{Nj} \big\} \quad \text{for }N,j\ge1.$$
 
 For $N$ fixed, $(\cG^j_N,\,j\ge1)$ is an exchangeable sequence of hierarchies on $[N]$. We define a random measure $\cL^N$ on the set of hierarchies on $[N]$ by specifying that for each such hierarchy $h$,
 \begin{equation}
  \cL^N\{h\} := \lim_{j\upto\infty} \frac{1}{j} \sum_{i\in [j]} \cf\left\{\cG^i_N = h\right\}.\label{eq:EIG_cnvgc}
 \end{equation}
 By de Finetti's theorem, this limit converges for each $h$. For each $N$, the resulting measure $\cL^N$ is a random probability measure. Moreover, de Finetti's theorem indicates that, given $\tail(\H_n)$, the sequence $(\cG^j_N,\,j\ge 1)$ is conditionally i.i.d.\ with law $\cL^N$.
 
 By the preceding argument, $\cL^N\{h\}$ can also be found by looking at asymptotic frequencies along a subsequence: 
 for each hierarchy $h$ on $[N]$,
 \begin{equation*}
 \begin{split}
  \cL^N\{h\} &= \lim_{k\upto\infty} \frac{1}{k} \sum_{i\in [k]} \cf\left\{\cG^{i(N+1)+1}_N = h\right\}
  		= \lim_{k\upto\infty} \frac{1}{k} \sum_{i\in [k]} \cf\left\{\Restrict{\cG^{iN+1}_{N+1}}{[N]} = h\right\}\\
  			&= \cL^{N+1}\{\text{hierarchies }h'\text{ on }[N+1]\colon \restrict{h'}{[N]} = h\} \quad \text{a.s.}.
 \end{split}
 \end{equation*}
 Thus, the family $(\cL^N,\,N\ge1)$ is a.s.\ projectively consistent. By the Daniell-Kolmogorov extension theorem, these laws extend to a random law $\cL$ on hierarchies on $\N$. Moreover, by essentially the same argument as for projective consistency, $\cL$ is almost surely an e.i.g.\ law: under $\cL^{kN}$ the hierarchies $(\cG^j_N,\,j\in [k])$ are i.i.d.\ with law $\cL^N$, and this holds for every $k\ge1$. We conclude that $\cL$ is an $\EIG$-valued regular conditional distribution for $(\H_n,\,n\ge1)$ given $\tail(\H_n)$.
 
 (iii) Suppose $\mathcal{M}$ is a random member of $\EIG$ and an r.c.d.\ for $(\H_n,\,n\ge1)$ over some $\sigma$-algebra. Let $\mathbf{Q}$ denote the distribution of $\mathcal{M}$ on $\EIG$, so
 $$\Pr\{(\H_n)\in \cdot\,\} = \EV[\mathcal{M}(\,\cdot\,)] = \int_{\EIG} \nu(\,\cdot\,)d\mathbf{Q}(\nu).$$
 For $j\ge1$ and $\nu\in\EIG$ deterministic, let $\nu^j$ denote the law of $\cG_j$ when $(\cG_n,\,n\ge1) \sim \nu$. In the notation of the proof of (ii), if $(\H_n,\,n\ge1)\sim\nu$ then $\cL^j = \nu^j$ a.s.. Thus, $\Pr(\cL^j = \mathcal{M}^j) = \int_{\EIG} 1 d\mathbf{Q}(\nu) = 1$. 
 Since laws of random hierarchies on $\N$ are specified by projections, $\mathcal{M} = \mathcal{L}$ a.s.. \qed



\section{Preliminaries for tree representations}\label{sect-preliminaries}

\subsection{Line-breaking construction of real trees}\label{sec:line_break}

\begin{definition}\label{def:l1}
 Let $\ell_1$ denote the Banach space of absolutely summable sequences of reals. Let $(\Be_j, j \geq 1)$ be the coordinate vectors in $\ell_1$, so that $\Be_1 = (1, 0, 0, \ldots)$, \mbox{$\Be_2  = (0,1,0,\ldots)$}, etc., and for $m \geq 1$ let $\pi_m$ denote the orthogonal projection onto span$\{\Be_1, \ldots, \Be_m\}$, and let $\pi_0$ send everything to $(0,0,\ldots)$, which we denote $0$.
 
 Let \textnormal{cl} denote the topological closure map on subsets of $\ell_1$.
\end{definition}

The following example of a real tree construction, due to Aldous \cite{MR1207226}, should provide adequate background for our purposes. For additional background, we refer the reader to \cite{MR2351587,MR2203728}.

\begin{example}[Line-breaking construction]\label{eg:line_break}
 Let $(L_n)$ be a sequence of positive numbers (not necessarily summable).  We define a family of real trees as follows: first, let $x_1 = 0\in \ell_1$ and let 
 \[
  \T_1 = x_1 + \Be_1[0,L_1] := \big\{(0, 0, \ldots) + \Be_1 z\colon 0 \leq z \leq L_1\big\}.
 \]
 Next, select a point $x_2$ from $\T_1$ and let
 \[
  \T_2= \T_1 \cup (x_2+\Be_2[0,L_2]) := \T_1 \cup \{x_2 + \Be_2 z\colon 0 \leq z \leq L_2\}.
 \]
 We continue recursively: supposing $\T_k$ has been defined, we select a point $x_{k+1}$ from $\T_k$ and set
 \[
  \T_{k+1} = \T_k \cup \big(x_{k+1}+\Be_{k+1}[0,L_{k+1}]\big),
 \]
 and let $\T = \cl(\bigcup_{n \geq 1} \T_n)$. The tree $\T_k$ is built up by ``gluing together'' $k$ line segments. Then $(\T,\ell_1,0)$ is a rooted real tree.
 
 To get a random real tree, simply randomize the construction above.  For example, let $(L_k)$ be the inter-arrival times of a Poisson process on $[0,\infty)$ of rate $t \, dt$, and for $k \geq 2$ select $x_{k}$ according to normalized length measure on $\T_k$.  The resulting tree is Aldous's Brownian continuum random tree.
\end{example}

All real trees constructed in this paper either result from the above construction or are non-branching trees on a single real interval. For $(\T,\ell_1,0)$ constructed in this way, for $x\in\T$, the segment from $0$ to $x$ is as follows.

\begin{definition}\label{def:special_path}
 Following Aldous \cite{MR1207226}, for $x\in \ell_1$ let $[[0,x]]_{sp}$ denote the path that proceeds from 0 to $x$ along successive directions. I.e.\ $[[0,x]]_{sp} := \textnormal{cl}[[0,x]]_{sp}^\circ$, where 
 \begin{equation}\label{special-path}
  [[0,x]]^\circ_{sp} := \bigcup_{m \geq 0} \{t\pi_m(x)+ (1-t)\pi_{m+1}(x)\colon t\in [0,1]\}.
 \end{equation}
\end{definition}

Observe that $[[0,x]]_{sp}$ differs from $[[0,x]]_{sp}^\circ$ only when $x=(x_1, x_2, \ldots)$ does not terminate in zeros, i.e when $x_j >0$ for infinitely many $j$, in which case $[[0,x]]_{sp} = [[0,x]]_{sp}^{\circ} \cup \{x\}$. We describe branch points and general segments, not just from $0$, in Definition \ref{def:tree_wedge}.

Recall Example \ref{eg:naive} and Figure \ref{fig:naive}. That figure exhibits three behaviors that may occur in exchangeable hierarchies on $\N$. From left to right, we call these iterative branching, broom-like explosion, and comb-like erosion. We now describe rooted, weighted real trees that correspond to these behaviors.

\begin{example}
 \begin{enumerate}[label=(\roman*), ref=(\roman*)]
  \item Consider $(\{0\},d,0,\delta_0)$, where $d$ is the metric $d(0,0) = 0$ and $\delta_0$ is a Dirac point mass. If $(\H_n,\,n\ge1)$ has law $\Theta(\{0\},d,0,\delta_0)$, then it a.s.\ equals the trivial hierarchy $\Xi(\N)$ where, for an arbitrary set $S$, the \emph{trivial hierarchy on $S$} is
  \begin{equation}\label{defn-xi}
   \Xi(S) := \{S\} \cup \{\{s\}\colon s\in S\} \cup \{\varnothing\}.
  \end{equation}
  This is an example of \emph{broom-like explosion}.
  
  \item Consider $([0,1],d,0,\Leb)$, where $d$ is the Euclidean metric and $\Leb$ is Lebesgue measure. The hierarchy obtained by sampling from this tree exhibits \emph{comb-like erosion}.
  
  \item Let $\mathcal{B}_{0} := \bigcup_{n\ge1}\{0,1\}^n$ and define $\phi\colon \mathcal{B}_0\to \N$ via $(b_i,\,i\in [n]) \mapsto \sum_{i=1}^n 2^{i-1}b_i$. Let $\mathcal{B}_{1} := \{0,1\}^{\N}$. For $b\in \mathcal{B}_1$ and $n\ge1$, let $b^{(n)}$ denote the truncation of $b$ to its first $n$ entries. We define $\eta\colon \mathcal{B}_{1} \to \ell_1$ via $\eta(b) := \sum_{n\ge1} 2^{-n}\Be_{\phi(b^{(n)})+1}$. Let $\T := \bigcup_{x\in\mathcal{B}_1} [[0,\eta(x)]]_{sp}$. Let $\psi\colon [0,1]\to\mathcal{B}_1$ denote the binary expansion map. Let $p$ denote the pushforward of Lebesgue measure via $\eta\circ\psi$. Then $(\T,\ell_1,0,p)$ is a rooted, weighted real tree, and the hierarchy obtained by sampling from it exhibits \emph{iterative branching}.
 \end{enumerate}
\end{example}

Here is an example in which two of these behaviors -- comb-like erosion and broom-like explosion -- are mixed together in a complex way.

\begin{example}[Fat Cantor weighted real tree]\label{eg:fat_Cantor}
 Let $A_0 := [0,1]$. Let $A_1 := A_0\setminus (3/8,5/8)$. We carry on recursively, as follows. For $n\geq 1$, $A_n$ comprises $2^n$ disjoint closed intervals of the same length. We form $A_{n+1}$ by removing an open interval of length $4^{-n-1}$ from the middle of each component of $A_n$. This sequence decreases to a fat Cantor set $A_{\infty} = \bigcap_{n\ge1}A_n$, also called a Smith-Volterra-Cantor set, with Lebesgue measure $1/2$; see \cite[p.\ 89]{MR1996162}.
 
 Let $B$ denote the set of maximal open intervals that comprise the complement of $A_{\infty}$. This is the set of intervals deleted in the course of the construction. 
 Let $d$ be Euclidean distance on $[0,1]$, $p_c$ the Lebesgue measure restricted to $A_{\infty}$, and $p_a := \sum_{(a,b)\in B}(b-a)\delta_a$. 
  Then $([0,1],d,0,p_a+p_c)$ is a rooted, weighted real tree.
\end{example}


\subsection{Most recent common ancestors and spinal compositions}\label{sec:MRCA}


\begin{definition}\label{MRCA-defn}
 If $\H$ is a hierarchy on a finite set $S$, then for $x,y\in S$, the \emph{most recent common ancestor (MRCA)} of $x$ and $y$ is
 \begin{equation} \label{defn-MRCA}
  (x \wedge y) := \bigcap_{G\in \H\colon x,y \in G} G.
 \end{equation}
 
 If $(\H_n,\,n\ge1)$ is a hierarchy on $\N$ then denote by $(i\wedge j)_n$ the MRCA of $i$ and $j$ in $\H_n$ if $n\ge i,j$, or the empty set otherwise. Denote by $(i\wedge j)$ the union $\bigcup_n (i\wedge j)_n$. 
 When discussing more than one hierarchy, e.g.\ $(\cG_n)$ and $(\H_n)$, we may write $(i \wedge j)_{\cG_k}$ or $(i \wedge j)_{\cG}$ to denote the MRCA of $i$ and $j$ in $\cG_k$ or in $(\cG_n,\,n\ge1)$, respectively.  
\end{definition}

Hierarchies on $\bN$ are fully specified by MRCAs.

\begin{proposition}\label{prop:MRCA}
 \begin{enumerate}[label=(\roman*), ref=(\roman*)]
  \item If $\H$ is a hierarchy on a finite set $S$ then \label{item:MRCA:cover}
  \[
   \H = \{(x \wedge y)\colon x, y \in S\}\cup \{\varnothing\}.
  \]
  \item If $(\H_n, n \geq 1)$ is a hierarchy on $\mathbb{N}$ then for every $n\geq \max\{i,j\}$\label{item:MRCA:intersect}
  \[
   (i \wedge j)_n = (i \bw j) \cap [n].
  \]
 \end{enumerate}
\end{proposition}

\begin{proof}
 (i) Note that $\{G\in \H\colon x\in G\}$ is totally ordered by inclusion, by part (a) of Definition \ref{defn-hierarchy}.  The smallest member of this class that contains $y$ is then $(x \wedge y)$.  This shows that
 \[
  \H \supseteq \{(x \wedge y)\colon  x, y \in S\}\cup \{\varnothing\}.
 \]
 To prove the reverse inclusion, fix non-empty $B \in \H$ and $x \in B$.  The class $\{(x\wedge y)\colon  y \in B\}$ is totally ordered by inclusion, with maximal element $(x \wedge y')$, say.  Then for all $z \in B$, $z \in (x \wedge z) \subseteq (x \wedge y')$, so $B \subseteq (x \wedge y')$.  On the other hand, $x,y' \in B$ and therefore $(x \wedge y')\subseteq B$.  This proves the reverse inclusion.  
 
 (ii) Since $\H_n = \H_{n+1}\big{|}_{[n]}$, for every $n \geq \max \{i,j\}$,
 \[
  [n] \cap (i\wedge j)_{n+1} = [n] \cap \bigcap_{G \in \H_{n+1}\colon  \{i,j\}\subseteq G} G = \bigcap_{G \in \H_{n}\colon  \{i,j\}\subseteq G}G = (i\wedge j)_{n}.
 \]
 By an inductive argument, $(i\wedge j)_n = (i\wedge j)_{N}\cap [n]$ for every $N\geq n$. The claim now follows from the definition of the MRCA $(i\wedge j)$ in $(\H_n)$.
\end{proof}


As a consequence of Definition \ref{defn-hierarchy}, for $i\in\N$ fixed, the MRCAs $\{(i\wedge j)\colon j\in\bN\}$ associated with a hierarchy $(\H_n,\,n\geq 1)$ are totally ordered by inclusion.

\begin{definition}
 A \emph{composition} of $\bN$ is a partition $\Pi$ of $\bN$ along with a total order $\preceq$ on the blocks. Let $(\H_n)$ be a hierarchy on $\bN$. For $i\in\bN$, we define
 \begin{equation}
 \begin{split}
  B_i(j) &:= (i\wedge j) \setminus \displaystyle\bigcup_{k\colon j\notin (i\wedge k)}(i\wedge k) \quad \text{for }j\in\bN\setminus\{i\}\text{ and}\\
  \Pi_i &:= \{ B_i(j)\colon j\in\bN\setminus\{i\}\}.
 \end{split}
 \end{equation}
 We say $B_i(j)\preceq_i B_i(k)$ if $(i\wedge j)\supseteq (i\wedge k)$. It follows from the definition of the MRCA that $(\Pi_i,\preceq_i)$ is totally ordered; we call this the \emph{$i^{th}$ spinal composition} with respect to $(\H_n)$.
\end{definition}

\begin{figure}[tbhp] 
   \centering
   \includegraphics[width=4in]{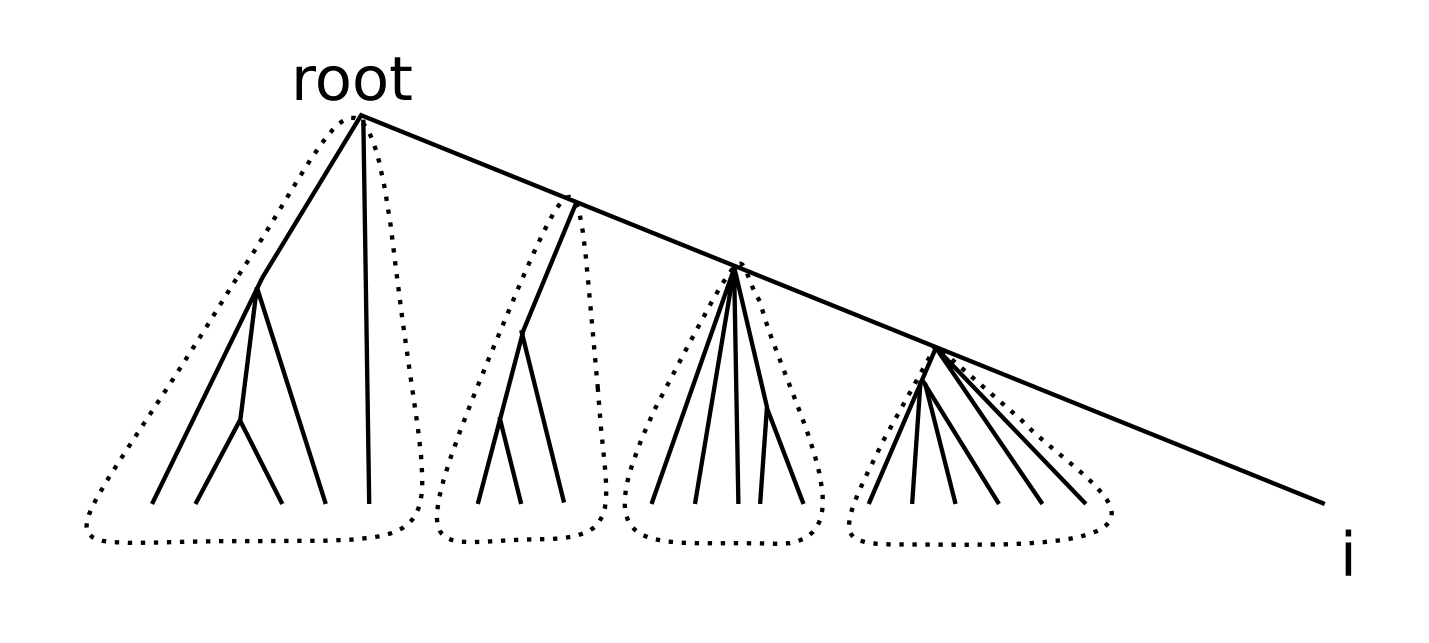} 
   \caption{The $i^{\rm{th}}$ spinal composition associated to a hierarchy is the partition of leaves of the hierarchy into blocks according to attachment point on the spinal path from root to leaf $i$, together with the following ordering on these blocks: block $s$ precedes block $t$ if the attachment point for block $s$ is nearer the root than the attachment point for block $t$.}
   \label{fig-spinal-comp}
\end{figure}

This can be described less formally in terms of the graph of $\H_n$ defined in Section \ref{sect-intro}. See Figure \ref{fig-spinal-comp}.
\begin{quote}
 For $i\in [n]$ draw the path from root to leaf $i$ in the graph of $\H_n$. The blocks of $\restrict{\Pi_i}{[n]}$ correspond to the subtrees sticking out from this path, and $B_i(j) \prec_i B_i(k)$ if the subtree corresponding to $B_i(j)$ is closer to the root. Now if we hold $i$ fixed and allow $n$ to increase without bound, these blocks grow and new blocks form, but order is preserved and blocks do not merge or fragment, resulting in a composition on $\bN\setminus \{i\}$ in the limit.
\end{quote}
%
If $(\H_n)$ is exchangeable then so are the associated spinal compositions. Exchangeable compositions have been studied previously and admit the following de Finetti-type description.

\begin{proposition}[a special case of \cite{MR1457625} Theorem 11, \cite{MR1104078} Theorem 5]\label{prop:spinal}
 Let $(\H_n)$ be an exchangeable hierarchy on $\bN$.
 \begin{enumerate}[label=(\roman*), ref=(\roman*)]
  \item For $i,j\in \bN$, the following limit exists almost surely:\label{item:spinal:def}
  \begin{equation}\label{eq:spinal_def}
   X^i_j := 1 - \lim_{n \to \infty} \frac{\#(i\wedge j)_n}{n}.
  \end{equation}
  We call these limits \emph{spinal variables}.
  \item For finite permutations $\sigma$ of $\bN$ (i.e.\ bijections with finitely many non-fixed points),\label{item:spinal:exch}
   \begin{equation}\label{eq:spinal_exch}
    \left(X^i_j;\ i,j\in\bN,\,i\neq j\right) \stackrel{d}{=} \left(X^{\sigma(i)}_{\sigma(j)};\ i,j\in\bN,\,i\neq j\right).
   \end{equation}
   In particular, for $i\in\bN$, the family $(X^i_j,\,j\in\bN\setminus\{i\})$ is exchangeable.
  \item For $i,j,k\in\bN$, the following events are almost surely equal:\label{item:spinal:equiv}
   \begin{equation}\label{eq:spinal_equiv}
    \{X^{i}_j \leq X^{i}_k\} = \{(i \bw k) \subseteq (i \bw j)\} = \{k \in (i\bw j)\}.
   \end{equation}
  \end{enumerate}
\end{proposition}

\begin{proof}[sketch]
 (i) The family $(\cf\{k\in (i\wedge j)\},\ k\in\bN\setminus\{i,j\})$ is exchangeable. Thus, the convergence follows from de Finetti's Theorem.
 
 (ii) This follows from the exchangeability of $(\H_n)$ and the definition of the spinal variables $(X^i_j)$.
 
 (iii) The only assertion in \eqref{eq:spinal_equiv} that doesn't follow easily from definitions is that if $k \notin (i \bw j)$ (or equivalently, if $(i\wedge k)\nsubseteq (i\wedge j)$) then $X^i_j > X^i_k$ almost surely. This can be deduced from applications of de Finetti's Theorem to the sequences
 \begin{gather*}
  \big(\cf\{(i\wedge k) \setminus (i\wedge j) = \{k\}\},\ k\in\bN\setminus\{i,j\}\big)\\
  \text{and} \quad \big(\cf\{l\in (i\wedge k)\setminus (i\wedge j)\},\ l\in\bN\setminus\{i,j,k\}\big).
 \end{gather*}
 In particular, the first sequence can't have more than a single `1,' or else $(\H_n)$ is not a hierarchy; and by de Finetti's theorem, it therefore a.s.\ contains no `1's. Consequently, for $k\notin (i\wedge j)$, the second sequence cannot be all zeroes, and thus must a.s.\ admit a positive limiting proportion of `1's. This limiting proportion equals $X^i_j - X^i_k$.
\end{proof}

We note that for all $i,j\in\N$ we have 
\begin{equation}
\begin{split}
 X^i_j = X^j_i \quad
 \text{and} \quad (i \wedge j) = \{m \in \bN\colon  X^i_m \geq X^i_j\}.\label{eq:spinal_MRCA}
\end{split}
\end{equation}
The former follows from the definition in \eqref{eq:spinal_def}. The latter restates part of \eqref{eq:spinal_equiv}.

\begin{corollary}\label{cor:all_big}
 Let $(\H_n)$ be an exchangeable hierarchy on $\bN$. For $i,j\in \bN$ distinct, $X^i_j < 1$ a.s.. Moreover, for each (non-random) infinite set $A\subseteq\bN$, the intersection $(i \wedge j) \cap A$ is a.s.\ infinite.
\end{corollary}

\begin{proof}
 Fix $i\neq j$. From \eqref{eq:spinal_equiv} we deduce that, since $j\notin (i\wedge i) = \{i\}$, we have $X^i_j < X^i_i = 1$, as desired. Thus, a positive asymptotic proportion of indices belong to $(i\wedge j)$. 
 The second claim then follows from de Finetti's theorem applied to $(\cf\{k\in (i\wedge j)\},\ k\in\bN\setminus\{i,j\})$.
\end{proof}

\begin{proposition}[Three-Rule for spinal variables]\label{prop:3rule}
 Let $(\H_n)$ be an exchangeable hierarchy on $\mathbb{N}$, with $(X^i_j)$ as in \eqref{eq:spinal_def}. For $i,j,k\in \bN$, if $X^i_k < X^j_k$ then $X^i_j = X^i_k$. Moreover, 
 there is a.s.\ some permutation $(i',j',k')$ of these indices for which
 \begin{equation}
  X^{i'}_{j'} = X^{i'}_{k'} \leq X^{j'}_{k'}.\label{eq:spinal_triple}
 \end{equation}
\end{proposition}

\begin{proof}
%
 It is easy to confirm that the two assertions are equivalent. Therefore, we prove only the former. Assume $X^{i}_k < X^j_{k}$. From \eqref{eq:spinal_equiv} we have $(j \wedge k) \subset (i \wedge k)$. This tells us: (a) $j\in (i\wedge k)$ and (b) $i\notin (j\wedge k)$. We apply \eqref{eq:spinal_equiv} to (b) twice, first to get $X^j_k > X^j_i$, and from there $k\in (i\wedge j)$. From this and (a) we conclude that $(i\wedge j) = (i\wedge k)$. Finally, by \eqref{eq:spinal_equiv}, $X^i_j=X^i_k$.
\end{proof}

\begin{figure}
 \centering
 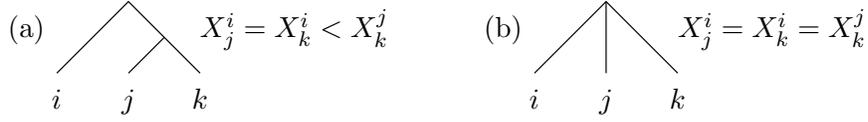
 \caption{Two distinct hierarchies on three elements correspond to cases of \eqref{eq:spinal_triple}.}\label{fig:spinal-triple}
\end{figure}

We call this proposition the \emph{Three-Rule} because it expresses in terms of spinal variables that, up to permutation of indices, there are only two possible hierarchies on a set of size three, as illustrated in Figure \ref{fig:spinal-triple}.

\section{Constructing a real tree from a hierarchy}\label{sec:construction}

Informally, in this section, we construct a random rooted, weighted real tree representation of an exchangeable hierarchy on $\N$, thereby proving Theorem \ref{thm:tree_constr}. Throughout this section, we consider only real trees constructed via line-breaking, as in Example \ref{eg:line_break}. In particular, these trees are rooted at $0 := (0,0,\ldots)$ and have the $\ell_1$ metric. Therefore, for brevity, we denote a weighted real tree simply $(\T,p)$, omitting the root and metric.

Let $(\H'_n, n \geq 1)$ be an exchangeable hierarchy on $\mathbb{N}$.  For reasons that will soon become clear, we prefer to work with a hierarchy on $\bZ$ rather than on $\bN$. Let $b: \bN\to \bZ$ denote the bijection that sends odd numbers to sequential non-positive numbers and evens to sequential positive numbers. For every $n \geq 1$ set 
\begin{equation}\label{eq:hier_on_Z}
 \H_n := \left\{\left\{b(k)\colon k \in (i \bw j)_{\H'}\right\} \cap [\pm n]\colon i,j \in \mathbb{N} \right\} \cup \Xi([\pm n]),
\end{equation}
where $[\pm n] := \{-n, \ldots, 0, \ldots, n\}$.  Then $\H_n$ is a hierarchy on $[\pm n]$ and $\Restrict{\H_{n+1}}{[\pm n]}  = \H_{n}$ for every $n \geq 1$.  We will need the notion of MRCAs in $\H_n$ and in $(\H_n,\,n\ge1)$; happily, Definition \ref{MRCA-defn} makes sense in the present context with obvious minimal changes, e.g.\ reading $[\pm n]$ for $[n]$.

Our main tool for constructing a real tree from $(\H_n)$ is the collection of $[0,1]$-valued spinal variables associated to spinal compositions:
\begin{equation}\label{eq:spinal_def_new}
 X^i_j := 1 - \lim_{n \to \infty}\frac{\#(i \bw j)_{\H_n}}{2n},   \qquad i,j \in \mathbb{Z}.
\end{equation}
Again, the results of Section \ref{sect-preliminaries} remain true here with obvious minimal changes. 
%
%
Our proof of Theorem \ref{thm:tree_constr} is organized as follows.
\begin{description}
 \item[\S\ref{sec:tree_constr}:] For $k<0$ we define real trees $\T_k$ and sequences $(t^k_j,\, j \geq 1)$ of random elements of $\T_k$, based upon $(\H_n)$. Informally, $\T_{-1}$ and $(t^{-1}_j)$ are a segment and samples corresponding to the $-1^{\rm{st}}$ spinal composition; $\T_{-2}$ has an additional branch, corresponding to the $-2^{\rm{nd}}$ spinal composition, splitting off from $\T_{-1}$ at a branch point corresponding to $(-1\wedge -2)_{\H}$; and so on. The $\T_k$ are shown to be the steps of a line-breaking construction building towards a limiting tree $\T$ bearing limiting samples $(t_j)$. We show that the $(t_j)$ are exchangeable, with a driving measure $p$.
 \item[\S\ref{sec:constr_eg}:] We give two illustrative examples of the construction.
 \item[\S\ref{sec:tree_MRCA}:] For distinct $u,v\in\bN$ we define a point $(t_u\wedge t_v)_{\ell}\in\T$ -- a branch point, except in degenerate cases -- that separates $t_u$ and $t_v$ from each other and the root. We show that the fringe subtree descending from this point corresponds to $(u\wedge v)_{\H}$.
 \item[\S\ref{sec:derived_hierarchy}:] Using the result of \S\ref{sec:tree_MRCA}, we prove that $\H_n\big|_{\bN}$ almost surely equals the hierarchy derived from $\T$ via the samples $(t_j)$. Since $\H_n\big|_{[n]}$ has the same distribution as $\H'_n$, this will complete our proof.
\end{description}


\subsection{Construction of trees and samples}\label{sec:tree_constr}


Recall the notation of Definitions \ref{def:l1} and \ref{def:special_path} for a standard basis $(\Be_n,\,n\ge1)$, projection maps $(\pi_n,\,n\ge1)$, and segments $[[0,x]]_{sp}$ in $\ell_1$.

\begin{definition}\label{def:constr_samples}
 For all $j \in \bZ$, set $t^{0}_j = 0$ and for every $k\le 0$,
 \begin{equation}\label{t-defn}
 \begin{split}
  t^{k-1}_j &:= t^k_j + \Be_{|k-1|}\left(X^{k-1}_j - \left\|t^k_j\right\|\right)_+,\\
  \T_k &:= \textnormal{cl}\left(\bigcup\nolimits_{j\geq 1} \left[\!\left[0, t^k_j\!\right]\right]_{sp} \right),
 \end{split}
 \end{equation}
 where $(a)_+ := \max\{a,0\}$. We treat $0$ as the root of each of the trees.
\end{definition}

We are mainly interested in the families $(t^k_j,\,j\geq 1)$. Definition \ref{def:constr_samples} can be described as follows: to define the samples $(t^{k-1}_j)$ for some $k\leq -1$, we select a subset of the $(t^k_j)$, possibly empty, and push these out in the $\Be_{|k-1|}$-direction, orthogonal to $\T_k$. For example, trivially, $\pi_{|i|}(t^k_j) = t^i_j$ for all $k<i<0$. Note that 
\begin{equation}
 \left\|t^k_j\right\| = \max_{i\in [k,-1]}X^i_j.\label{eq:mag_tnk}
\end{equation}

We will now show that all of the samples that are pushed out in passing from $t^k_j$ to $t^{k-1}_j$ are selected from the same spot, $t^k_{k-1}$, on $\T_k$. Thus, $\T_{k-1}$ is derived by adding at most one branch to $\T_k$. We call this the \emph{Line-Breaking Property} because it shows that our sequential construction of the $\T_k$ fits into the framework of the line-breaking construction of Example \ref{eg:line_break}.

\begin{lemma}[The Line-Breaking Property of $\T$]\label{lem:its_a_tree}
 For $k\leq -1$ and $j\geq 1$, if $t^{k-1}_j\neq t^k_j$ then $t^k_j = t^k_{k-1}$. Moreover, regardless of whether $t^{k-1}_j=t^k_j$,
 \begin{equation}
  \left(X^{k-1}_j - \left\|t^k_j\right\|\right)_+ = \left(X^{k-1}_j - \left\|t^k_{k-1}\right\|\right)_+.\label{eq:new_branch_X}
 \end{equation}
\end{lemma}

\begin{proof}
 Note that $t^{k-1}_j\neq t^k_j$ if and only if $X^{k-1}_j > \|t^k_j\|$. By \eqref{eq:mag_tnk}, this means $X^{k-1}_j > X^i_j$ for $i\in [k,-1]$. By the Three-Rule \eqref{eq:spinal_triple}, this implies $X^i_j = X^i_{k-1}$ for each such $i$; so by definition, $t^k_j = t^k_{k-1}$. This proves the first assertion, as well as \eqref{eq:new_branch_X} in the case $X^{k-1}_j > \|t^k_j\|$.
 
 Now suppose $X^{k-1}_j \leq \|t^k_j\|$. Then by \eqref{eq:mag_tnk} there is some $i\in [k,-1]$ for which $X^{k-1}_j \leq X^i_j$. By the Three-Rule \eqref{eq:spinal_triple}, this implies $X^{k-1}_j \leq X^i_{k-1}$, so again by \eqref{eq:mag_tnk}, $X^{k-1}_j \leq \|t^k_{k-1}\|$.
\end{proof}

\begin{proposition}\label{prop:tree_limits}
 $\T_{-1} \subseteq \T_{-2} \subseteq \cdots$. Moreover, the limits 
 \begin{equation}
  t_j := \lim_{k \to -\infty} t^k_j \quad (j \geq 1)
 \end{equation}
 exist and are members of $\T := \textnormal{cl}(\T^{\circ})$, where
 \begin{equation*}
  \T^{\circ} := \bigcup_{k<0} \T_k.
 \end{equation*}
 Finally, $\T$ is a random real tree.
\end{proposition}

\begin{proof}
 By Lemma \ref{lem:its_a_tree} and Definition \ref{def:special_path} of the segments $[[0,x]]_{sp}$, for every $j\ge1$ and $k<0$, $[[0,t^k_j]]_{sp}\subseteq [[0,t^{k-1}_j]]_{sp}$. Thus, $\T_k\subseteq\T_{k-1}$. By definition, the spinal variables $(X^i_j)$ take values in $[0,1]$ almost surely, so $\|t^k_j\| \leq 1$. Since, $\pi_{|k|}(t^{k-1}_j) = t^k_j$, this gives the desired convergence results. Finally, by Lemma \ref{lem:its_a_tree}, $\T$ is a real tree resulting from a line-breaking construction, as in Example \ref{eg:line_break}.
\end{proof}

\begin{figure}[hbtp] 
   \centering
   \includegraphics[width=3in]{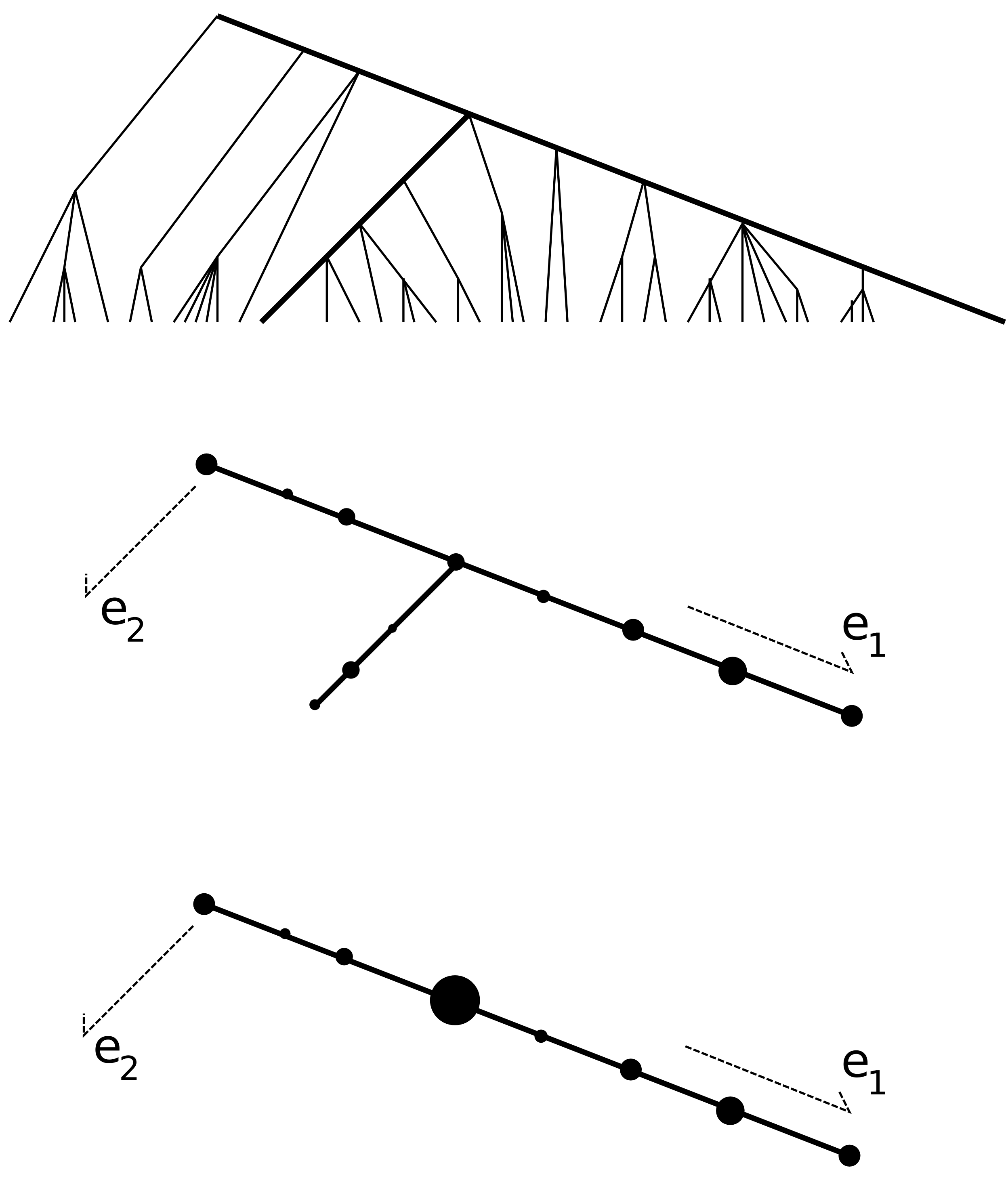}
   \caption{At top is shown the graph of $\H_n$ with leaf labels erased.  The bold paths are the spinal paths to leaves $-1$ and $-2$, respectively.  In the middle, $(\T_{-2}, p_{-2})$ is shown.  The arrows indicate the $\ell_1$ basis directions, and atoms of $p_{-2}$ are represented by black circles or {\em beads} on $\T_2$, with circle size corresponding to atom size.  At bottom is shown $(\T_{-1},p_{-1})$.  Note that $(T_{-2}, p_{-2})$ is derived from $(\T_{-1}, p_{-1})$ by ``crushing'' a bead on $\T_{-1}$ into fragments and stringing the crushed bead fragments out in the $\Be_{2}$ direction.}
   \label{fig:bead_crush}
\end{figure}

This construction is an example of \emph{bead crushing}; see \cite{MR2561439}. This is illustrated in Figure \ref{fig:bead_crush}.

\begin{proposition}\label{prop:driving}
 The family $(t_j,\,j\geq 1)$ is exchangeable and has a driving measure $p$. Likewise, for every $k<0$, the family $(t^k_j, j\geq 1)$ is exchangeable and has a driving measure $p_k$.
\end{proposition}

\begin{proof}
 The exchangeability of the $(t^k_j,\,j\geq 1)$ follows from that of the vectors $\big((X^i_j,\; i\in [k,-1]),\; j\geq 1\big)$, noted in \eqref{eq:spinal_exch}. The $(t_j)$ are exchangeable since they are the limits of the $(t^k_j)$. By de Finetti's Theorem, these sequences have driving measures on $\ell_1$.
\end{proof}

\subsection{Examples of the construction}\label{sec:constr_eg}

The following examples may help the reader visualize the preceding construction. 

\begin{example}
 Let $(U_n, n \in \mathbb{Z})$ be a family of i.i.d.\ Uniform[0,1] random variables and let 
 \[
  \H_n := \big\{\{j \in [\pm n]\colon U_j \geq x\}\colon 0 \leq x \leq 1\big\} \cup  \Xi([\pm n]).
 \]
 Following the construction in Section \ref{sec:tree_constr}, $\T_{-1} = \Be_1[0,U_{-1}]$, each $t^{-1}_j = \min\{U_{-1},U_j\}$, and $p_{-1}$ is length measure on $\T_{-1}$ plus an atom of mass $1-U_{-1}$ at $\Be_1U_{-1}$.  Now let $k_{1}=-1$ and define a sequence $(k_{m}, m \geq 1)$ recursively by $k_{m+1}:=\max\{i<0\colon U_{i} > U_{k_m}\}$.  Then 
 \begin{gather*}
  (\T_{-1},p_{-1}) = \cdots = (\T_{k_2 +1}, p_{k_2+1})\\
  \text{and} \quad \T_{k_2} = \T_{-1} \cup \left(\Be_1U_{-1} + \Be_{|k_2|}[0,U_{k_2}-U_{k_1}]\right).
 \end{gather*}
 I.e.\ $\T_{k_2}$ is an isometric embedding of $[0,U_{k_{2}}]$ in $\ell_1$, with a {\em kink} or bend at the image of $U_{-1}$ in $\ell_1$. The measure $p_{k_2}$ is the length measure on $\T_{k_2}$ plus an atom of size $1-U_{k_2}$ at the far end, $\Be_{k_1}U_{k_1} + \Be_{k_2}(U_{k_2}-U_{k_1})$.  In general, $\T_{k_m}$ is an isometric embedding of $[0,U_{k_m}]$ into $\ell_1$ with $|k_m| - 1$ kinks, and $p_k$ is length measure on $\T_{k_{m}}$ plus an atom of size $1-U_{k_m}$ at the end of $\T_{k_m}$. The limit tree $\T$ is an isometric copy of $[0,1]$, embedded in $\ell_1$, and $p$ is length measure on $\T$.  This tree has only one leaf, which has $p$-measure 0.
\end{example}

\begin{example}
 Let $(B_n,\,n\in\mathbb{Z})$ be i.i.d.\ Bernoulli trials with success rate $1/2$. Let
 \[
  \H_n := \big\{ \{j\in [\pm n]\colon B_j = b\}\colon b=0,1\big\} \cup \Xi([\pm n]).
 \]
 Then for every $i\neq j$, if $B_i=B_j=b$ then $(i\wedge j)$ is the set of all $k\in\mathbb{Z}$ for which $B_k = b$, and if $B_i\neq B_j$ then $(i\wedge j) = \mathbb{Z}$. Thus, $X^i_j = \frac12\cf\{B_i = B_j\}$. Therefore, $\T_{-1} = \Be_1[0,1/2]$ and $t^{-1}_j = \Be_1\frac12\cf\{B_i = B_j\}$ for each $j\neq -1$. Consequently, $p_{-1} = \frac12\delta_0 + \frac12\delta_{\Be_1/2}$.
 
 Now, let $k = \max\{i<-1\colon B_k\neq B_{-1}\}$. Then for $-1\ge i>k$ we have $(\T_i,p_i) = (\T_{-1},p_{-1})$, and
 $\T_k = \Be_1\left[0,\frac12\right] \cup \Be_{|k|}\left[0,\frac12\right]$. 
 For $i\notin [k,-1]$, if $B_i = B_1$ then $t^k_i = \frac12\Be_1$, and otherwise $t^k_i = \frac12\Be_{|k|}$. So $p_k = \frac12\delta_{\Be_1/2} + \frac12\delta_{\Be_{|k|}/2}$. Finally, $(\T_i,p_i) = (\T_k,p_k)$ for all $i<k$.
\end{example}

\subsection{MRCAs correspond to points in the tree}\label{sec:tree_MRCA}

\begin{definition}\label{def:tree_wedge}
 For $x,y\in \ell_1$ with all non-negative coordinates,
 $$[[0,x]]_{sp}\cap [[0,y]]_{sp} = [[0,z]]_{sp}$$
 for some $z\in \ell_1$, possibly equal to zero. We define 
 \begin{equation}\label{eq:tree_wedge_def}
  (x\wedge y)_\ell := z, \qquad [[x,y]]_{sp} := \big([[0,x]]_{sp} \cup [[0,y]]_{sp}\setminus [[0,z]]_{sp}\big) \cup \{z\}.
 \end{equation}
\end{definition}

\begin{figure}[htbp]
 \centering
 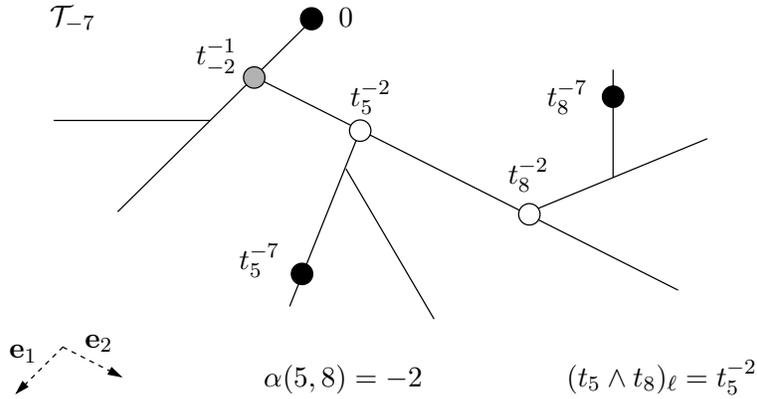
 \caption{Important vertices for Lemma \ref{lem:BP_position}.}
 \label{fig:branch_pt}
\end{figure}

In the notation of \cite{MR1207226}, $(x\wedge y)_\ell$ is instead called $b(x,y)$. For points $x,y\in\T$, if $x\notin [[0,y]]_{sp}$ and $y\notin [[0,x]]_{sp}$ then $(x\wedge y)_{\ell}$ is the branch point of $\T$ that separates $x$ and $y$ from $0$ and each other; i.e.
\begin{equation*}
 [[0,x]]_{sp} \cap [[0,y]]_{sp} \cap [[x,y]]_{sp} = \{(x\wedge y)_{\ell}\}.
\end{equation*}
And the path $[[x,y]]_{sp}$ is in fact the segment in $\T$ between $x$ and $y$, as in Definition \ref{def:real_tree}. In the example in Figure \ref{fig:branch_pt}, $(t_5\wedge t_8)_{\ell} = t^{-2}_5$.

\begin{proposition}\label{prop:MRCA_equals_fringe}
 For distinct $u,v\in\bN$ and $\fringe{\cdot}{\T}$ as in \eqref{fringe-subtree-defn},
 \begin{equation}
  (u\wedge v)_{\H}\cap\N = \{j\in\N\colon t_j\in \fringe{(t_u\wedge t_v)_{\ell}}{\T}\} \quad \text{a.s.}.
  \label{eq:MRCA_equals_fringe}
 \end{equation}
\end{proposition}

Towards this result, we require a pair of lemmas.

\begin{lemma}\label{lem:alpha_in_MRCA}
 Let $(x_i,\,i\geq 1)$ and $(y_i,\,i\geq 1)$ denote the $\ell_1$ coordinates of $t_u$ and $t_v$ respectively, for some distinct $u,v\in\bN$. Then
 \begin{equation}\label{eq:last_+_coord}
  \alpha(u,v) := -\max(\{1\}\cup \{i\geq 2\colon \min(x_i,y_i) > 0\})
 \end{equation}
 is a.s.\ finite, with $\alpha(u,v)\in (u\wedge v)_{\H}$.
\end{lemma}

\begin{proof}
 Let $\alpha := \alpha(u,v)$ and let
 \begin{equation*}
  \beta := \max\left((u\wedge v)_{\H}\cap \bZ_-\right).
 \end{equation*}
 By Corollary \ref{cor:all_big}, $(u\wedge v)_{\H}\cap \bZ_-$ is not empty, so $\beta$ is finite, well-defined, and belongs to $(u\wedge v)_{\H}$. We will prove $\alpha=\beta$ a.s..
 
 
 By \eqref{eq:spinal_equiv}, the event $\{k\in (u\wedge v)_{\H}\}$ is equivalent to both of $\{X^u_k \geq X^u_v\}$ and $\{X^v_k\geq X^v_u\}$, so
 \begin{equation}\label{eq:second_form_beta}
  \beta = \max\{k < 0\colon X^k_u \geq X^u_v\} = \max\{k < 0\colon X^k_v\geq X^u_v\}.
 \end{equation}
 If $\beta = -1$ then $\beta \geq \alpha$ trivially. Otherwise, if $\beta < -1$, then the maximality of $\beta$ implies
 \[
  X^{\beta}_u > \max\{X^{\beta+1}_u,\ldots,X^{-1}_u\} = \|t^{\beta+1}_u\|,
 \]
 and correspondingly for $X^\beta_v$. Thus, by \eqref{t-defn} and \eqref{eq:mag_tnk}, both $t_u$ and $t_v$ have a positive $|\beta|^{\rm{th}}$ coordinates. By definition of $\alpha$, this means $\beta\geq \alpha$.
 
 To prove $\alpha = \beta$ it will suffice to show that for $k<\beta$, if $t_u$ has a positive $|k|^{\rm{th}}$ coordinate then $t_v$ does not. So suppose such a $k$. Then 
 \begin{equation*}
 \begin{split}
  X^k_u > \|t^{k+1}_u\| \ge X^{\beta}_u \ge X^u_v = X^k_v,
 \end{split}
 \end{equation*}
 where the first inequality follows from \eqref{eq:mag_tnk}, the second from \eqref{t-defn}, the third from \eqref{eq:second_form_beta}, and the last then follows by the Three-Rule \eqref{eq:spinal_triple}. With another appeal to \eqref{eq:second_form_beta}, this gives $X^k_v = X^u_v \leq X^{\beta}_v$. Thus, by \eqref{t-defn}, the $|k|^{\rm{th}}$ coordinate of $t_v$ is 0.
\end{proof}

\begin{lemma}\label{lem:BP_position}
 Fix $u,v\in\bN$ distinct and let $\alpha := \alpha(u,v)$, as in \eqref{eq:last_+_coord}. Then
 \begin{align}
  (t_u\wedge t_v)_{\ell} &= t^{\alpha+1}_{\alpha} + \Be_{\alpha}\left(\min\left\{X^{\alpha}_u,X^{\alpha}_v\right\} - \left\|t^{\alpha+1}_{\alpha}\right\|\right)\label{eq:tree_MRCA}\\
  	&= t^{\alpha+1}_{\alpha} + \Be_{\alpha}\left(X^u_v - \left\|t^{\alpha+1}_{\alpha}\right\|\right) \quad \text{a.s.}\label{eq:tree_MRCA_2}\\[.1cm]
  \text{and} \quad 
  X^u_v &= \|(t_u\wedge t_v)_{\ell}\| \ge \left\|t^{\alpha+1}_\alpha\right\| \quad \text{a.s.}.\label{eq:BP_position}
 \end{align}
 Moreover, the final inequality is sharp only when $X^u_v = 0$.
\end{lemma}


\begin{proof}
%
 \eqref{eq:tree_MRCA}: In the trivial case in which there is no $k<0$ for which $t_u$ and $t_v$ both have positive $|k|^{\text{th}}$ coordinates, we get $\alpha=-1$ and \eqref{eq:tree_MRCA} holds with both sides equal to $0$. Now suppose there is some such $k$. By the Line-Breaking Property, Lemma \ref{lem:its_a_tree}, since $t_u$ and $t_v$ have positive $|\alpha|^{\rm{th}}$ coordinates, $t^{\alpha+1}_u = t^{\alpha+1}_v = t^{\alpha+1}_{\alpha}$. This means that $t^{\alpha}_u$ and $t^{\alpha}_v$ differ only in their $|\alpha|^{\text{th}}$ coordinates, if at all. Thus, \eqref{eq:tree_MRCA} follows from \eqref{eq:new_branch_X} and the definitions of $(\,\cdot\wedge\cdot\,)_{\ell}$ and the $(t^k_j)$.
 
 \eqref{eq:tree_MRCA_2}: By Lemma \ref{lem:alpha_in_MRCA}, $\alpha\in (u\wedge v)_{\H}$. Thus, by \eqref{eq:spinal_MRCA}, $\min\{X^{\alpha}_u,X^{\alpha}_v\}\ge X^u_v$. Applying the Three-Rule, \eqref{eq:spinal_triple}, we get $\min\{X^{\alpha}_u,X^{\alpha}_v\} = X^u_v$. Substituting this into \eqref{eq:tree_MRCA} gives \eqref{eq:tree_MRCA_2}.
 
 \eqref{eq:BP_position}: This follows immediately from \eqref{eq:tree_MRCA_2} and the definition of $\alpha$. 
\end{proof}

\begin{proof}[Proof of Proposition \ref{prop:MRCA_equals_fringe}]
 Let $\alpha := \alpha(u,v)$. First, suppose $(t_u\wedge t_v)_{\ell}\in [[0,t_j]]_{sp}$, and we will show $j\in (u\wedge v)_{\H}$. 
 From Definition \ref{def:special_path} of $[[0,t_j]]_{sp}$, the coordinates of $t_j - (t_u\wedge t_v)_{\ell}$ are all non-negative. Thus, by \eqref{t-defn} and \eqref{eq:tree_MRCA_2}, $X^{\alpha}_j \geq \min\{X^{\alpha}_u,X^{\alpha}_v\}$. 
 By \eqref{eq:spinal_MRCA}, $j$ must belong to at least one of $(u\wedge\alpha)_{\H}$ or $(v\wedge \alpha)_{\H}$. Since we also have $\alpha\in (u\wedge v)_{\H}$ from Lemma \ref{lem:alpha_in_MRCA}, by Definition \ref{MRCA-defn} of MRCAs, $(u\wedge\alpha)_{\H}\cup (v\wedge \alpha)_{\H}\subseteq (u\wedge v)_{\H}$. Thus, $j\in (u\wedge v)_{\H}$, as desired.
 
 
 Now suppose instead $j\in (u\wedge v)_{\H}$ and we will show $(t_u\wedge t_v)_{\ell}\in [[0,t_j]]_{sp}$. If $X^u_v = 0$ then, by \eqref{eq:BP_position}, $(t_u\wedge t_v)_{\ell} = 0$. Since $0\in [[0,t_j]]_{sp}$, this case is trivial, so we can assume $X^u_v>0$. 
 By Lemma \ref{lem:alpha_in_MRCA}, $\alpha\in (u\wedge v)_{\H}$, and so $(\alpha\wedge j)_{\H} \subseteq (u\wedge v)_{\H}$. 
 Therefore, $X^\alpha_j \geq X^u_v$. By the sharpness condition for \eqref{eq:BP_position}, $X^\alpha_j > \|t^{\alpha+1}_{\alpha}\|$. By the Line-Breaking Property, Lemma \ref{lem:its_a_tree}, $t^{\alpha+1}_j = t^{\alpha+1}_{\alpha}$. Finally, by \eqref{eq:tree_MRCA_2},
 $$(t_u\wedge t_v)_{\ell} = t^{\alpha+1}_{\alpha} + (X^u_v - \|t^{\alpha+1}_{\alpha}\|)\Be_{\alpha}\quad \text{whereas} \quad t^\alpha_j = t^{\alpha+1}_{\alpha} + (X^\alpha_j - \|t^{\alpha+1}_{\alpha}\|)\Be_{\alpha}.$$
 We conclude that $(t_u\wedge t_v)_{\ell} \in [[0,t^{\alpha}_j]]_{sp}\subseteq [[0,t_j]]_{sp}$, as desired.
\end{proof}

\begin{corollary}\label{cor:p_diffuse_boundary}
 The measure $p$ is a.s.\ diffuse on $\partial\T := \T\setminus\T^\circ$.
\end{corollary}

\begin{proof}
 Take $u,v\in\bN$ distinct. If $t_u = t_v$ then $t_u = (t_u\wedge t_v)_{\ell}$. By Lemma \ref{lem:alpha_in_MRCA}, $\alpha(u,v) > -\infty$. By \eqref{eq:tree_MRCA_2}, this point $(t_u\wedge t_v)_{\ell}$ has only finitely many non-zero coordinates; i.e.\ $(t_u\wedge t_v)_{\ell}\notin\partial\T$. The claim follows since $p$ is the driving measure of the $(t_j)$. 
\end{proof}

\subsection{The derived hierarchy}\label{sec:derived_hierarchy}

For $n \geq 1$ and $x\in \T_k$ we define
\begin{equation}
\begin{split}
 I_n(x) &:= \{j \in [n]\colon t_j\in \fringe{x}{\T}\} = \{j \in [n]\colon x\in [[0,t_j]]_{sp}\},\label{eq:In_def_1}\\
 \cI_n &:= \{I_n(y)\colon y\in \T\} \cup \Xi([n]).
\end{split}
\end{equation}
This $\cI_n$ is the hierarchy derived from $\T$ via the samples $(t_1,\ldots, t_n)$.

\begin{proposition}\label{prop:In_Gn}
 $\cI_n= \H_n\big|_{[n]}$ almost surely, for every $n\in\bN$. 
\end{proposition}

\begin{proof}
 Let $\cG_n := \H_n\big|_{[n]}$. Propositions \ref{prop:MRCA}\ref{item:MRCA:cover} and \ref{prop:MRCA_equals_fringe} have the immediate consequence $\cG_n\subseteq \cI_n$. It remains to prove the reverse inclusion: we must show $I_n(x)\in \cG_n$ for each $x\in\cT$. Consider the subtree $\T^n := \bigcup_{j=1}^n [[0,t_j]]_{sp} \subseteq \T$. If $x\notin\T^n$ then $I_n(x) = \varnothing$, so we need only consider $x\in\T^n$.
 
 We call the $(t_u,\,u\in [n])$ and the $((t_u\wedge t_v)_{\ell},\ u,v\in [n])$ \emph{marked points} in $\T^n$. There are finitely many marked points. By definition, the leaves and branch points of $\T^n$ are all marked. Thus, if $x\in \T^n$ is not a marked point then it lies along some branch, say in the $|k|^{\rm{th}}$ coordinate direction. Then there is some nearest marked point $z$ beyond $x$ with $x\in [[0,z]]_{sp}$. I.e.\ this point lies along the same branch as $x$, with no samples or branch points in between, so $I_n(x) = I_n(z)$. Now, it suffices to show $I_n(z) \in \cG_n$ when $z$ is a marked point.
 
 Case 1: $z = (t_u\wedge t_v)_{\ell}$ for some distinct $u,v\in [n]$. Then $I_n(z) = (u\wedge v)_{\cG_n} \in \cG_n$ by Proposition \ref{prop:MRCA_equals_fringe}.
 
 Case 2: $z = t_u$ for some $u\in [n]$. If $I_n(t_u) = \{u\}$ then we're done: $\{u\}\in \cG_n$. Otherwise, take $v\in I_n(t_u)\setminus\{u\}$. Then by definition, $t_u\in [[0,t_v]]_{sp}$, and so $t_u = (t_u\wedge t_v)_{\ell}$, thus reducing the problem to Case 1. 
\end{proof}

\begin{proof}[Proof of Theorem \ref{thm:tree_constr}]
 By Theorem \ref{thm:EIG}, our original hierarchy $(\H'_n,\,n\ge1)$ on $\N$ admits a r.c.d.\ $\cL$ given $\tail(\H'_n)$, where $\cL$ is a random member of $\EIG$. By exchangeability, if $(\H'_n)$ has law $\cL$ then $(\H_n|_{[n]},\,n\ge1)$ has law $\cL$ as well. Thus, by Theorem \ref{thm:EIG}\ref{item:EIG:unique}, $\cL$ is a.s.\ unique with the property of being a random member of $\EIG$ that is a r.c.d.\ for $(\H_n|_{[n]},\,n\ge1)$.
 
 By Proposition \ref{prop:driving}, $p$ is a driving measure for the samples $(t_j,\,j\ge1)$. Thus, $\Theta(\T,0,\ell_1,p)$ is a r.c.d.\ for $(\cI_n,\,n\ge1)$, and by Proposition \ref{prop:In_Gn}, for $(\H_n|_{[n]},\,n\ge1)$ as well. Since $\Theta(\T,0,\ell_1,p)$ is a random member of $\EIG$ as well, it equals $\cL$ a.s.. Thus, it is a r.c.d.\ for $(\H'_n,\,n\ge1)$ given $\tail(\H'_n)$. 
\end{proof}



\section{Proof of Theorem \ref{thm:interval_Kingman}}\label{sec:interval_results}

 Let $(\H_n)$ be an exchangeable hierarchy on $\bN$ and let $(\T,\ell_1,0,p)$ denote the corresponding tree constructed in Section \ref{sec:tree_constr}, so that $\Theta(\T,\ell_1,0,p)$ is a r.c.d.\ for $(\H_n)$. 
 For $k \in \mathbb{N}$ let $\pi_k$ be as in Definition \ref{def:l1}. We consider $\ell_1$ to be totally ordered lexicographically:
 \begin{quote}
  $(x_i,\, i\geq 1) \prec (y_i,\, i\geq 1)$ if there is some $k\ge 0$ such that $\pi_k(x) = \pi_k(y)$ and $x_{k+1} < y_{k+1}$.
 \end{quote}
 
 We will map the set of fringe subtrees of $\T$ to an interval hierarchy by pulling back via a process $\xi\colon [0,1]\to \T$ that explores $\T$ in lexicographic order. This is similar to a ``wall follower algorithm'' or depth-first search, always turning to explore branches added later in the iterative construction of $\T$ before jumping back to continue along a main branch. To that end, let $\T_k := \pi_k(\T)$ and $p_k := \pi_k(p)$. For $x\in \T$, the set $\{y\in \T\colon y\preceq x\}$ is Borel. Thus, for $k\ge 1$ we may define
 \begin{equation*}
 \begin{split}
  D_k(x) &:= p_k\{y\in \T_k\colon y\preceq x\} \quad \text{for }x\in\T_k\\
  \text{and} \quad \xi_k(u) &:= \sup\{x\in \T_k\colon D_k(x) < u\} \quad \text{for }u\in [0,1].
 \end{split}
 \end{equation*}
 To visualize this, we refer back to the ``bead-crushing'' view of the construction of $(\T,p)$, illustrated in Figure \ref{fig:bead_crush} and described in Lemma \ref{lem:its_a_tree}. For each $k\ge1$, $p_{k+1}$ is formed from $p_k$ by reducing or eliminating the mass at one atom at a point $x\in\T_k$, and replacing that mass to lie along a new branch $B = \T_{k+1}\setminus \T_k$ rooted at $x$. 
 %
 Suppose, for example, that $p_3$ has an atom of mass $1/2$ at a point $x\in\T_3$, and in step $4$ this atom is partially crushed so that $p_4$ still has an atom at $x$, but with mass only $1/4$, with another $1/4$ mass in $p_4$ spread out along a new, linear branch $B = \T_4\setminus\T_3$ rooted at $x$. The lexicographic order on $\T_4$ is formed by taking the lexicographic order on $\T_3$ and inserting $B$ so that it follows immediately after the point $x$: for $x\prec y\in \T_3$ and $z\in B$ we have $x\prec z\prec y$. Thus, $\xi_3$ will pause for time $1/2$ at $x$ -- say, for example, on the time interval $(1/4,3/4]$ -- whereas $\xi_4$ will only pause at $x$ on $(1/4,1/2]$, then spending time $(1/2,3/4]$ exploring the branch $B$.
  
 It follows from this construction that for $j<k$ we have $D_j(x)\ge D_k(x)$ for $x\in \T_j$ and $\pi_j(\xi_k(u)) = \xi_j(u)$ for $u\in [0,1]$. Thus, we can define
 \begin{align}
  D(x) &:= \lim_{k\to\infty} D_k(x) = p\{y\in \T\colon y\preceq x\} \quad \text{for }x\in\bigcup\nolimits_k\T_k,\\
  \text{and} \quad \xi(u) &:= \lim_{k\to\infty}\xi_k(u) = \sup\left\{x\in \bigcup\nolimits_k\T_k\colon D(x) < u\right\} \quad \text{for }u\in [0,1].\notag
 \end{align}
 
 By Corollary \ref{cor:p_diffuse_boundary}, $p$ is diffuse on $\partial\T := \T\setminus\bigcup_k\T_k$. For $x\in \partial\T$, the fringe subtree $\fringe{x}{\T} = \{x\}$. For $x\in\bigcup_k\T_k$, let $K(x)$ denote the least index $k$ for which $x\in \T_k$. Each fringe subtree equals the intersection of $\T$ with an interval in the lexicographic order, so $D$ maps it to a real interval:
 \begin{equation}\label{eq:lexico_fringe}
 \begin{split}
  \fringe{x}{\T} &= \{y\in\T\colon x\preceq y \preceq x+\Be_{K(x)}\},\\
  \text{so} \quad D(\fringe{x}{\T}) &= \big[D(x),D(x)+p(\fringe{x}{\T})\big],\\
  \text{and thus} \quad \fringe{x}{\T} &\supset \xi\big(D(x),D(x)+p(\fringe{x}{\T})\big].
 \end{split}
 \end{equation}
 The fringe subtrees $\fringe{x}{\T}$ with $x\in \bigcup_k\T_k$ comprise a $\pi$-system that generates the Borel $\sigma$-algebra on $\T$. Via Dynkin's $\pi$-$\lambda$ Theorem, we conclude that $p$ is the push-forward of Lebesgue measure under $\xi$.
 
 In light of \eqref{eq:lexico_fringe} and for consistency with Definition \ref{def:derived_by_sampling}\ref{item:int_hier}, let
 \begin{equation}
  \mathscr{H}:= \left\{\big[D(x),D(x)+p(\fringe{x}{\T})\big)\colon x \in \T\right\} \cup \Xi([0,1)).
 \end{equation}
 Let $(U_j,\,j\ge1)$ be an i.i.d.\ sequence of Uniform$[0,1)$ random variables independent of $\T$, and let $t_j = \xi(U_j)$ for each $j$, so the $(t_j,\,j\ge1)$ have driving measure $p$. Let $I_n(x)$ be as in \eqref{eq:In_def_1}. For $n \geq 1$ let 
 \begin{equation*}
  \H_n' := \{\{j \in [n]\colon U_j \in B\}\colon B \in \mathscr{H}\},
  \qquad \H''_n := \left\{ I_n(x)\colon x\in\T\right\} \cup \Xi([n]).
 \end{equation*}
 It is easily seen that $\H_n'' = \H_n'$ for every $n$, almost surely. Moreover, $(\H'_n,\,n\ge1)$ has r.c.d.\ $\Theta(\scH,\Leb)$ given $(\T,p)$, while $(\H_n'',\,n\ge1)$ has r.c.d.\ $\Theta(\T,\ell_1,0,p)$. Thus, $\Theta(\scH,\Leb) = \Theta(\T,\ell_1,0,p)$ a.s., and so it is a r.c.d.\ for $(\H_n,\,n\ge1)$ given $\tail(\H_n)$, as desired. \qed

\section{Complements}\label{sect-complements}

\subsection{Exchangeable hierarchy probability functions}\label{sect-ehpf}

Recall that the {\em graph} of a hierarchy $\H_n$ on $[n]$ is a rooted tree $\tT_n$ with $n$ leaves, and with certain other properties noted in the introduction, where each leaf bears a distinct label in $[n]$.  We define the {\em shape} of such a tree $\tT_n$ to be its orbit under the action of the symmetric group, the action of a permutation $\sigma$ being to relabel leaf $i$ by $\sigma(i)$ for every $i \in [n]$.  We use lower-case bold face $\mathtt{s}_n$ to denote the shape of $\tT_n$, which can be regarded as a function of $\tT_n$ or alternatively of $\H_n$:  $\mathtt{s}_n=\mathtt{s}(\tT_n) = \mathtt{s}(\H_n)$.  Obviously, the shape of $\tT_n$ can be identified with the unlabeled tree derived by erasing the labels on the leaves of $\tT_n$, but this observation is not far from a tautology, as unlabeled graphs are often defined as orbits under such actions of the symmetric group. 

We write $\tts_n \nearrow \tts_{n+1}$ if it is possible to remove a leaf of $\mathtt{s}_{n+1}$ and thereby obtain $\mathtt{s}_n$.  In this context, {\em removing a leaf} means (i) erasing the leaf and the edge of $\mathtt{s}_{n+1}$ which had that leaf as an endpoint and then (ii) if removing said leaf results in the leaf's parent now having only one remaining child -- i.e.\ degree two -- then contracting the edge between parent and child. Equivalently, $\tts_n \nearrow \tts_{n+1}$ if there is a hierarchy $\H_{n+1}$ on $[n+1]$ for which $\mathtt{s}_n = \mathtt{s}(\H_{n+1}|_{[n]})$ and $\mathtt{s}_{n+1} = \mathtt{s}(\H_{n+1})$.

Let $\mathtt{S}:= \{\mathtt{s}(\H_n)\colon \H_n$ a hierarchy on $[n]$ for some $n \geq 1\}$.  If $(\H_n)$ is an exchangeable hierarchy on $\mathbb{N}$, then there is a function $h\colon \mathtt{S} \mapsto [0,1]$ for which for every fixed hierarchy ${H}$ on $[n]$,
\begin{align}
 \Pr(\H_n = H) &= h(\mathtt{s}(H)),\label{symmetry-h}\\
 h\big(\mathtt{s}(\Xi(\{1\}))\big) &= 1,\label{normalization-h}\\
 \text{and} \quad h(\mathtt{s}_n) &= \sum_{\mathtt{s}_{n+1}\colon \mathtt{s}_n \nearrow \mathtt{s}_{n+1}} h(\mathtt{s}_{n+1}).\label{consistency-h}
\end{align}
Here, $\Xi(\{1\})$ is the trivial hierarchy on $\{1\}$. 
 Formula (\ref{symmetry-h}) asserts that $\Pr(\H_n = H)$ only depends on the shape of the graph of $H$. 
 We call $h$ the \emph{exchangeable hierarchy probability function} (EHPF) representation of the law of $(\H_n,\, n\ge1)$. This is analogous to the  representation of exchangeable partitions by EPPFs, discussed in \cite{MR1337249} and \cite[Chapters 2 \& 3]{MR2245368}, or that of exchangeable compositions via consistent composition structures \cite{MR1457625}.

Under the relation $\nearrow$, the space of rooted tree shapes becomes a graded graph or lattice, in the sense of \cite{MR1984868}. 
In this framework, the problem of characterizing exchangeable hierarchies on $\N$ (solved by Theorem \ref{thm:tree_constr}) is equivalent to characterizing the class of bounded, positive harmonic functions on this lattice. See \cite[Chapter 0]{MR1984868} for much more on this topic.

\begin{proposition}\label{prop:EHPF}
Suppose that $h\colon \mathtt{S}\mapsto [0,1]$ satisfies \eqref{consistency-h} for all $n \geq 1 $ and \eqref{normalization-h}. Then there is an exchangeable random hierarchy $(\H_n,\,n\ge1)$ on $\mathbb{N}$ for which \eqref{symmetry-h} holds for every fixed hierarchy $H$ on $[n]$, for all $n \geq 1$.
\end{proposition}
\begin{proof}
 Let $\H_1 = \Xi([1])$, and assuming that $\H_1, \ldots, \H_n$ have been defined, conditionally given $\H_n = H$, select $\H_{n+1}$ from the set 
 \[\left\{\text{hierachies }H'\text{ on }[n+1]\colon \Restrict{H'}{[n]} = H\right\},\]
 selecting $H'$ with probability $h(s(H'))/h(s(H))$.  
\end{proof}

Now suppose $(\H_n,\, n\ge 1)$ is exchangeable, with EHPF $h$, and consider sequences of disjoint, finite subsets $A_1,\ldots,A_k\subset [N]$, for some $N>k>1$, and corresponding hierarchies $H_1,\ldots,H_k$, with each $H_i$ being a hierarchy on $A_i$. By definition of the EHPF,
$$\Pr\left( \forall i\in [k],\ \restrict{\H_N}{A_i} = H_i\right) = \sum_{\text{hierarchies }H\text{ on }[N]\colon \forall j\in [k],\ \restrict{H}{A_j}=H_j}h(\mathtt{s}(H)).$$
If $(\H_n)$ is independently generated, this implies
\begin{equation}\label{eq:IG_EHPF}
 \sum_{\text{hierarchies }H\text{ on }[N]\colon \forall j\in [k],\ \restrict{H}{A_j}=H_j}h(\mathtt{s}(H)) = \prod_{i=1}^k h(\mathtt{s}(H_i)).
\end{equation}
Conversely, any EHPF that satisfies \eqref{eq:IG_EHPF} for every pair of sequences $(A_1,\ldots,A_k)$ and $(H_1,\ldots,H_k)$ describes an independently generated hierarchy. In light of this, Theorem \ref{thm:EIG} and Proposition \ref{prop:EHPF} have the following corollary.

\begin{corollary}
 All EHPFs, i.e.\ all functions $h$ that satisfy \eqref{consistency-h} for all $n \geq 1 $ and \eqref{normalization-h}, can be expressed uniquely as convex combinations of EHPFs that satisfy \eqref{eq:IG_EHPF} for all $(A_1,\ldots,A_k)$ and $(H_1,\ldots,H_k)$ as above.
\end{corollary}

\subsection{Open problems}\label{questions}

Recall the Na\"{\i}ve conjecture stated in the introduction, that the three behaviors appearing in Figure \ref{fig:naive} -- iterative branching, broom-like explosion, and comb-like erosion -- are in some sense the only behaviors that arise in exchangeable hierarchies. In light of Theorem \ref{thm:tree_constr}, we offer a formal interpretation of this.

\begin{conjecture}\label{conj:formal_naive}
 Let $(\H_n)$ be an exchangeable hierarchy on $\bN$, with $(\T,0,\ell_1,p)$ the random rooted, weighted real tree constructed from $(\H_n)$ as in Section \ref{sec:tree_constr}. Then $p$ can a.s.\ be uniquely decomposed into a sum of three components, $p = p_a+p_s+p_l$, where $p_a$ is purely atomic, $p_s$ is the restriction of length measure to a subset of the skeleton of $\T$, and $p_l$ is a diffuse measure supported on the leaves of $\T$.
\end{conjecture}

The purely atomic measure $p_a$ corresponds to explosions in the hierarchy, $p_s$ corresponds to erosion, and $p_l$ may be said to describe a part of the hierarchy that ``survives'' these explosions and erosion, remaining in large blocks that dwindle incrementally, via splitting. 


In addition to this conjecture, we propose two open-ended problems.\medskip

\noindent\emph{Problem 1.} 
%
For rooted, weighted real trees $(\T_1,d_1,r_1,p_1),\ (\T_2,d_2,r_2,p_2)$, 
is there a simple way to see whether $\Theta(\T_1,d_1,r_1,p_1) =$ $\Theta(\T_2,d_2,r_2,p_2)$? Speaking loosely, it is possible to {\em prune away} tree branches of $\T_1$ that carry no $p_1$-mass, and also stretch segments of $\T_1$ arbitrarily, and not change $\Theta(\T_1,d_1,r_1,p_1)$.  Purely topological considerations are not quite enough to settle this question: suppose that $\T_1$ is the tree [0,1] rooted at 0 and $p_1$ is Lebesgue measure on [0,1], and suppose that $\T_2$ is is the half line $[0,\infty)$ rooted at 0, and $p_2$ is the Exponential(1) distribution on $\T_2$.  Then $\T_1$ and $\T_2$ are not homeomorphic, but $\Theta(\T_1,d,0,p_1) = \Theta(\T_2,d,0,p_2)$.\medskip



\noindent\emph{Problem 2.} 
According to a result in \cite{MR1337249}, if $p$ is an EPPF then there is a sequence $(P_1, P_2, \ldots)$ of nonnegative random variables with $\sum_iP_i \leq 1$, such that for all $\lambda \in \bigcup_{k \geq 1} \mathbb{N}^k$,
\begin{equation}\label{pitman-moments}
p(\lambda_1, \ldots, \lambda_k) = \mathbb{E}\left[ \left( \prod_{i=1}^k P_i^{\lambda_i-1} \right) \left( \prod_{i=2}^k (1-P_1- \ldots - P_{i-1}) \right) \right] .
\end{equation}
The EPPF $p$ determines the joint law of this sequence $(P_i)$ uniquely, and conversely. 
There is also the following consequence of Theorem \ref{thm:partn_Kingman}. Let $m_\lambda$ denote the monomial symmetric polynomial 
\begin{equation}
m_{\lambda}(x):=\sum_{\sigma}x^{\lambda_1}_{\sigma(1)}\ldots x^{\lambda_k}_{\sigma(k)}
\end{equation}
for $\lambda = (\lambda_1, \ldots, \lambda_k)$, where the sum is taken over all injective functions $\sigma:[k] \mapsto \mathbb{N}$. Suppose $\Pi$ is an exchangeable random partition of $\mathbb{N}$ for which the sequence ($P_i$) of (\ref{pitman-moments}) satisfies
 $\Pr\big( \sum\nolimits_j P_j = 1\big) = 1$, 
i.e.\ $\Pi$ almost surely does not contain any singleton blocks. Then there is a measure $\mu$ on the Kingman simplex $\nabla:=\{(x_i,\,i\ge1)\colon x_1 \geq x_2 \geq \ldots,\,\sum_{i} x_i = 1\}$ such that
\begin{equation}\label{symmetric-moments}
 p(\lambda) = \int_{\nabla}m_{\lambda}(x) \mu(dx)  \qquad \text{for }\lambda\in\bigcup_{k \geq 1} \mathbb{N}^k.
\end{equation}
In particular, $\mu$ is the law of the rearrangement of the sequence $(P_i,\,i\ge1)$ in non-increasing order. For proofs of the preceding material, see \cite{MR1337249} and \cite[Chapters 2 \& 3]{MR2245368}.  

Does there exist an exchangeable hierarchies analogue to \eqref{symmetric-moments} that is reasonably comprehensive, e.g.\ that goes beyond describing hierarchies formed by well-ordered recursive splitting?


\vskip .4cm

We thank David Aldous, Steve Evans, and Matthias Winkel for helpful discussion. We also thank reviewers for their careful readings, both of the first version of the paper in 2011, as well as the current version in 2017.

\bibliographystyle{plain}
\bibliography{hierarchiesbib}

\end{document}